\documentclass[english,a4paper,11pt]{article}
\usepackage{amsmath,amssymb}
\usepackage{graphicx,color,url}
\usepackage{babel}
 
\newcommand{\RR}{\mathbb{R}}

\newcommand{\QQ}{\mathbb{Q}}

\newtheorem{Theorem}{Theorem}
\newtheorem{Proposition}[Theorem]{Proposition}
\newtheorem{Definition}{Definition}

\begin{document}
\title{Finite Orbits in Multivalued Maps and Bernoulli Convolutions} 
\author{Christoph Bandt\footnote{This work was supported by Deutsche Forschungsgemeinschaft grant Ba 1332/11-1.}
}
\maketitle\vspace{-6ex}

\begin{abstract}
Bernoulli convolutions are certain measures on the unit interval depending on a parameter $\beta$ between 1 and 2. In spite of their simple definition, they are not yet well understood. We study their two-dimensional density which exists by a theorem of Solomyak.   To each Bernoulli convolution, there is an interval $D$ called the overlap region, and a map which assigns two values to each point of $D$ and one value to all other points of $[0,1].$  There are two types of finite orbits of these multivalued maps which correspond to zeros and potential singularities of the density, respectively.

Orbits which do not meet $D$ belong to an ordinary map called $\beta$-transformation  and exist for all $\beta>1.6182.$ They were studied by Erd¨os, J´oo, Komornik, Sidorov, de Vries and others as points with unique addresses, and by Jordan, Shmerkin and Solomyak as points with maximal local dimension. In the two-dimensional view, these orbits form address curves related to the Milnor-Thurston itineraries in one-dimensional dynamics. The curves depend smoothly on the parameter and represent quantiles of all corresponding Bernoulli convolutions.  

Finite orbits which intersect $D$ have a network-like structure and can exist only at Perron parameters $\beta.$ Their points are intersections of extended address curves, and can have finite or countable number of addresses, as found by Sidorov. For an uncountable number of parameters, the central point $\frac12$ has only two addresses. The intersection of periodic address curves can lead to singularities of the measures. We give examples which are not Pisot or Salem parameters.
\end{abstract}
 
\section{Introduction}
This paper studies two related subjects: the occurence of network-like orbits in multivalued maps, and the parametric family of Bernoulli convolutions, BCs for short.  The orbit of a point $x$ in a dynamical system can be visualized as a path of arrows leading from $f^k(x)$ to $f^{k+1}(x)$ for $k=0,1,2,...$ The orbit becomes finite if $x$ fulfils the equation $f^{k+p}(x)=f^k(x)$ for some $k$ and $p,$ so that the path turns into a cycle. 
For a multivalued map, the orbit of $x$ is represented by a branching tree. Such orbit can only be finite if all branches lead back to lower levels of the tree, which is usually expressed by several equations, determining not only the point $x$ but also the mapping $f,$ at least to some extent. Nevertheless, network-like orbits appear in simple multivalued maps of an interval, in particular those which define BCs. 

\begin{figure}[h]
\begin{center}
\includegraphics[width=0.75\textwidth]{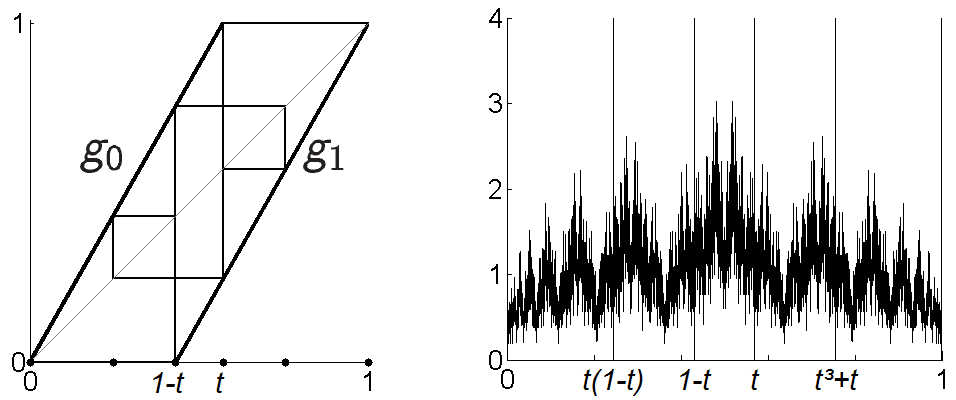}
\end{center}
\caption{The mappings $g_0(x)=\beta x, \, g_1(x)=\beta x +1-\beta$ and the corresponding Bernoulli convolution for $\beta \approx 1.7549,$ the Pisot number with minimal polynomial $\beta^3 -2\beta^2+\beta -1.$ The finite orbit of
$t=1/\beta$ is indicated here, together with the generated partition. }
\label{bc}
\end{figure}

Bernoulli convolutions are the simplest examples of self-similar measures with overlaps. They have been studied as examples in real analysis since the 1930s. Given a number $\beta$ between 1 and 2, let $t=1/\beta$ and consider the two linear functions 
\begin{equation} g_0:[0,t]\to [0,1], \ g_0(x)=\beta x\quad\mbox{ and }\quad g_1:[1-t,1]\to [0,1], \ g_1(x)=\beta x+1-\beta \, ,\label{expa}\end{equation} 
as indicated in Figure \ref{bc}.  The BC with parameter $\beta$ is the unique probability measure $\nu$ on $[0,1]$ which fulfils

\begin{equation} 
\nu (A)=\frac12 \nu(g_0(A))+\frac12 \nu(g_1(A))\quad\mbox{ for all Borel sets } A\subset [0,1] .\label{bern}\end{equation}

When $A$ is a subset of $[0,1-t]$ or $[t,1],$ only one term will appear on the right of \eqref{bern} since $g_1(A)$ or $g_0(A)$ is empty. The map $G(x)=\{ g_0(x), g_1(x)\}$ is multivalued only for $x\in [1-t,t]=D,$ the so-called overlap region. Basic facts on  BCs can be found in \cite{So4,PSS,ban} and in Section 6. Definitions of Pisot and Garsia numbers are given in Section 4.

Since Erd\"os \cite{E} proved in 1939 that  BCs for Pisot numbers $\beta$ are singular measures, much work was done to determine those $\beta$ for which $\nu$ is absolutely continuous.  There are detailed studies of  BCs of Pisot numbers \cite{AF,F4,FS,KL,HHM} and other algebraic numbers \cite{G,FW,F11,Va}. Singularity of $\nu$ is still only known for the countable set of Pisot numbers, and regularity of specific Bernoulli convolutions could be shown only for the countable set of Garsia numbers \cite{G}, not even for any particular rational number $\beta .$ However, Solomyak \cite{So} proved in 1995 that for Lebesgue almost all $\beta$ the measure $\nu$ has a density function, even an $L^2$ density function as shown by Peres and Solomyak \cite{PS}. 
In other words, if we take a random number $\beta$ between 1 and 2, then $\nu$ will have a density function with probability one. See the surveys \cite{PSS,So4} for more information on related work. Recently, Shmerkin \cite{Sh14} used new tools of Hochman to improve Solomyak's result further: the set of $\beta$ leading to singular  BCs has Hausdorff dimension zero. In 2015, Varju \cite{Va} proved that for all algebraic numbers $\beta\in [2-\varepsilon , 2]$ fulfilling a technical condition, the measure $\nu$ is absolutely continuous. However, $\varepsilon$ was not specified, it could be as small as $10^{-10^{10}}.$ More recent work by Breuillard and Varju \cite{BV} shows that the dimension of $\nu_t$ is 1 for transcendental numbers $t$ like ln 2, $e^{-\frac12},$ and $\pi/4,$ while Hare and Sidorov \cite{HS} proved that the dimension is 1 also for new algebraic parameters.

This paper follows a different path. We study $\nu$ for the interesting parameter region $\beta\in [1.618, 1.99]$ and include structural properties of $\nu$ obtained by  Erd\"os, J\'oo, Komornik, Sidorov, de Vries and others \cite{dVK,EJK,GS,JSS,KL,Si3,Si9} in connection with $\beta$-expansions. 
We consider all $\nu_t$ with $t=1/\beta \in (\frac12 , 1)$ together as a two-dimensional `super-measure'. Solomyak's work implies that this measure is given by an $L^2$ density function of two variables. The problem to classify BCs into `singular' and `regular' now turns into the study of singularities of an $L_2$ function on a rectangle. (The following statement is slightly stronger than Solomyak's result in \cite{So}, it follows from \cite{PS}, and in an obvious way from Shmerkin's \cite{Sh17} recent result that almost all BCs are $L^q$ for $q>1.$)

\begin{Theorem}[Two-dimensional density of Bernoulli convolutions \cite{So,PS}] \label{2D}\hfill\\
There is an $L^2$ function $\Phi: [\frac12,1]\times[0,1]\to [0,\infty)$ such that for Lebesgue almost all parameters $t=1/\beta\in  [\frac12,1],$ the density of the Bernoulli convolution $\nu_t$ is the function $\Phi(t,x), x\in [0,1].$
\end{Theorem}

\begin{figure}[h]
\includegraphics[width=0.999\textwidth]{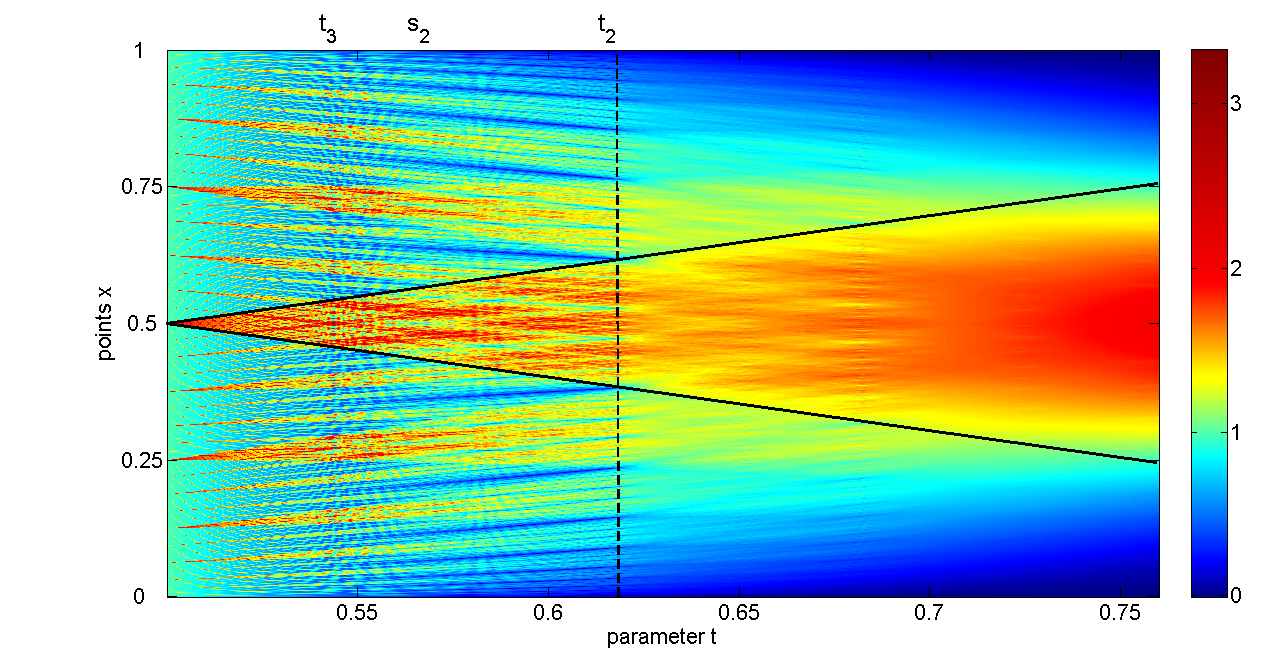}
\caption{The function $\Phi$ for $0.5\le t\le 0.76.$ Bernoulli convolutions for 4000 values $t=1/\beta$
were approximated by histograms, and visualized as vertical sections over $t$ with color code indicated on the right. The measure of Figure \ref{bc} appears at $t=s_2\approx 0.57.$ The overlap region $\bf D$ is between the black lines. We study the parameters $t<t_2\approx 0.618,$ left of the dashed vertical line.}
\label{B5}
\end{figure}

Figure \ref{B5} shows an approximation of $\Phi$ for $0.5\le t\le 0.76.$  The intricate structure in this picture will be studied in Sections 6-10  where magnifications of Figure \ref{B5} are shown. For each $t,$ the average value is 1. All values are between 0 and 3.3 as the color bar on the right indicates.  The small values of $\Phi ,$ drawn in dark blue, as well as the large values, drawn in red, follow certain patterns. 
The purpose of our paper is to clarify this structure, with the goal of better understanding the singularities of certain Bernoulli convolutions, and how the measures $\nu_t$ change with the parameter $t.$ 

In the middle of Figure \ref{B5} we see the triangular overlap region ${\bf D}=\{ (t,x)|\ t\le x\le 1-t\}$ which contains large values of $\Phi .$ In Section 6 we explain the smaller overlap regions which form `horns' with tip on the left axis. Next, there is a kind of phase transition at the Fibonacci parameter $t_2\approx 0.618.$ On the left of $t_2,$ dark blue curves separate the horns. Henceforth, we focus on this case $t<t_2$ where we have apparent structure.
The dark curves outside ${\bf D}$ become more numerous when $t$ goes down to $\frac12,$ and they seem to form a Cantor structure.  These dark curves consist of points with unique address in the fractal construction, which is discussed in Section 2 below.

We define $A_t$ as the set of points with unique address for the parameter $t.$ Here, $A_t$ can be viewed as
the intersection of the dark curves with the vertical line $x=t.$ These sets were studied by many authors, in connection with $\beta$-expansion.  Around 1990,
Erd\"os and co-authors started to consider those $\beta=1/t$ for which the equation $\sum_{k=1}^\infty c_kt^k=1$ has a unique solution $c_k\in\{ 0,1\}.$ They were called univoque numbers, and characterized in \cite{EJK} and \cite{AF}. In our Figure \ref{B5} the univoque parameters are those $t$ where the dark curves enter the overlap triangle ${\bf D},$ cf.  Proposition \ref{univoque}.  Furthermore, it was proved that $A_t$ is infinite for $t<t_2,$ see \cite{EJK} and its references, and Daroczy and Katai \cite{DaK,DaKa}. Komornik and Loreti \cite{KL} found the constant $\beta_{KL}\approx 1.7872$ at which $A_t$ becomes uncountable. Allouche and Cosnard \cite{AC} proved that this number is transcendental. Glendinning and Sidorov \cite{GS} formally defined the set $A_t,$ proved that it contains accumulation points for $t\le s_2\approx 0.57$ and is a Cantor set for $\frac12<t<t_{KL}.$ Further work includes \cite{ACS,AF,JSS,Kal,KoLi,Si2,Si3,Si9} and the comprehensive paper of de Vries and Komornik \cite{dVK} on the topology of $A_t$ and of the set of univoque numbers. 
 
In Section 7 and 8, we provide a self-contained approach to all these results. While the quoted authors used $\beta$-expansions, we define binary itineraries and address curves, similar to the technique introduced by Milnor and Thurston \cite{MT} for the study of unimodal maps. In our view, $A_t$ is the set on which $G(x)= \{ g_0(x), g_1(x)\}$ acts as an ordinary map. Our basic Proposition \ref{conjlemma} says that the action of $G$ on $A_t$ is conjugate to the doubling map $g(x)=2x \mod 1$ on a corresponding subset of $[0,1],$ for all $t<t_2.$ 

While previous authors studied Cantor sets $A_t\subset [0,1]$ for fixed $t,$ our focus is on address curves $x_b(t)$ for binary sequences $b$ in the two-dimensional setting. These curves form a smooth element in an otherwise chaotic scenario. In Theorem \ref{quant} we show that each address curve $x_b(t)$ describes the $b$-quantile of all BCs in its domain. Here $b$ is considered as binary number in $[0,1].$ All ordinary periodic and preperiodic points of the map $G$ are represented by address curves, and the parameters where such curves enter ${\bf D}$ are landmarks of $\Phi .$ 

Although $\Phi$ describes just a measure on a rectangle, we consider it as a representation of a parametric family of dynamical systems, like the Mandelbrot set, or the bifurcation diagram of the family of real quadratic maps  \cite{CE,PR}. Our approach with binary itineraries underlines the tight connection between these three objects which are all based on the doubling map $g(x)=2x \mod 1.$ Allouche, Clarke and Sidorov \cite{ACS} verified the Sharkovskii order for periodic points in $A_t.$ In Section 8 we show that unimodal maps and the sets $A_t$ are in a one-to-one correspondence provided by binary itineraries, which allows to transfer other results between the two theories. 

Thus the dynamics of $G$ is well understood on the sets $A_t$ where it is an ordinary map. The difficulty comes when we go inside the overlap region where $G$ is really multivalued. However, we can extend the address curves into ${\bf D}$ and consider their intersections. If at least one of the addresses is nonperiodic and the $G$-orbit of the intersection point $y$ does not return to the overlap interval, then $y$ has either two or a countable number of addresses (Propositions \ref{intersec1} and \ref{intersecaleph}). This implies that the local dimension of $\nu_t$ at $y$ assumes the same maximum as at the points with unique addresses. Such points were studied by Sidorov \cite{Si9}, see also \cite{Baker,BS}. In Section 9, we discuss this case, solve a problem of  \cite{Si9} and prove the following surprising fact (Theorem \ref{central}): for an uncountable number of parameters, the central point $\frac12$ has only two addresses.

\begin{figure}[h!]
\begin{center}
\includegraphics[width=0.75\textwidth]{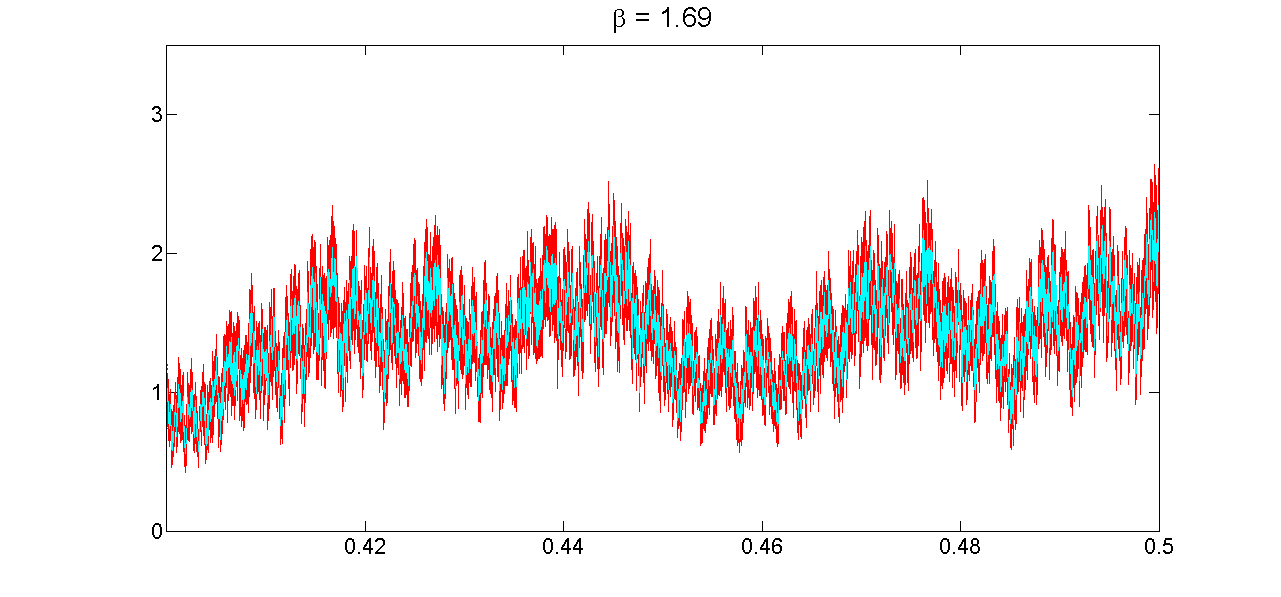}\\
\includegraphics[width=0.75\textwidth]{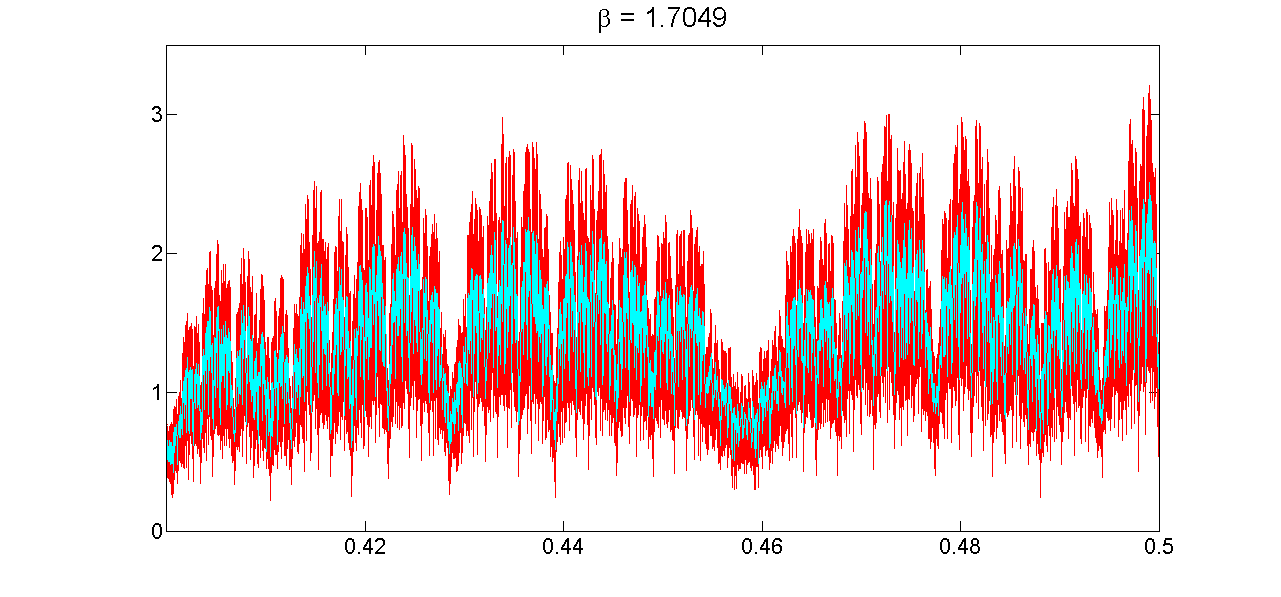}\\
\includegraphics[width=0.75\textwidth]{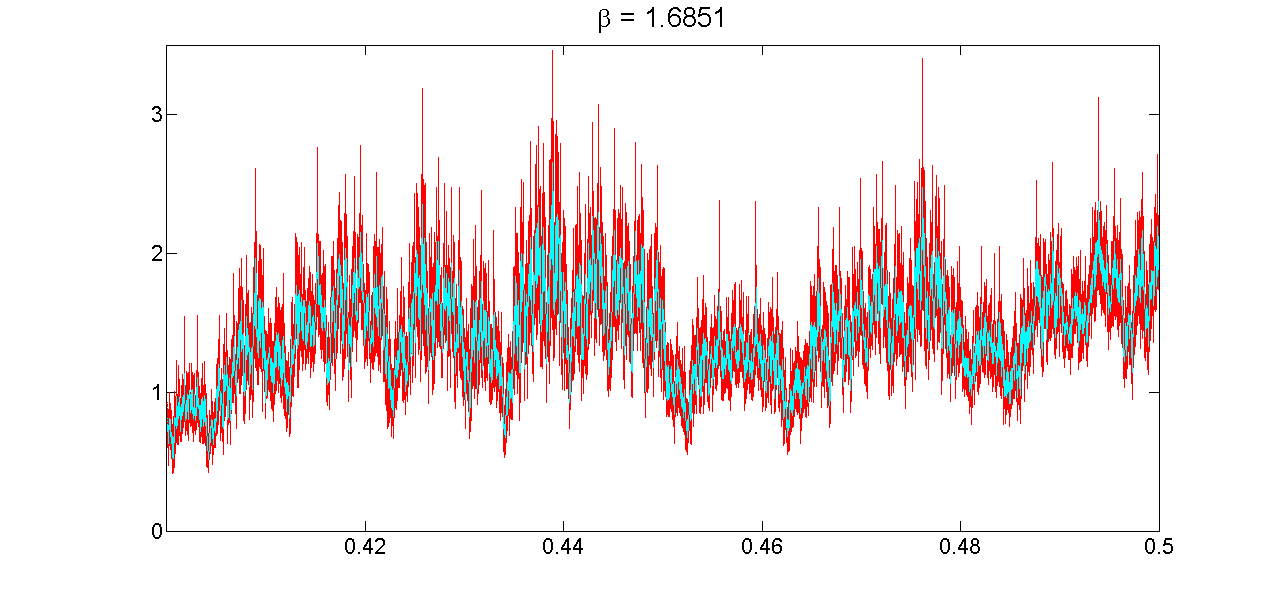}
\end{center}
\caption{An impression of Bernoulli convolutions for parameters $\beta=1.69$ (presumably a density function, even continuous), $\beta\approx 1.7049$ (Pisot parameter, singular measure) and $\beta=1.6851$ (Perron parameter, see Section 10). The Matlab figures include 5000 blue and 250000 red equidistant points in $[0.4, 0.5],$ adjacent values are connected by lines.}
\label{drei}
\end{figure}

The intersection of address curves is most interesting when both addresses are periodic. This case, discussed in Section 10, is the only one which can lead to singularities of $\Phi .$  Here the $G$-orbit returns to $D$ through at least two different cycles, and we can estimate the local dimension of $d=d_y(\nu)$ at the intersection point $y.$ In the supercritical case where $d>2t$ there is a singularity, which means that $\nu$ cannot have a bounded or continuous  density function (Theorem \ref{intersecper}). This includes a theorem of Feng and Wang \cite{FW}.
We provide new  examples of parameters with such weak singularities which are neither Pisot nor Salem. One of them is shown in Figure \ref{drei}: the Perron number $\beta\approx 1.6851$ has higher peaks than nearby Pisot parameters but probably the corresponding BC has an $L_1$ density. 

There are many `subcritical' intersections of periodic address curves where we expect a bounded density of $\nu$ and a nontrivial multifractal spectrum, as found by Feng \cite{F11} and Jordan, Shmerkin and Solomyak \cite{JSS}. Moreover, Theorem \ref{nonex} says that intersections of periodic or preperiodic address curves can happen only for parameters $\beta$ which are weak Perron numbers.  Thurston \cite{Th} pointed out the importance of weak Perron numbers in one-dimensional dynamics.\smallskip

The illustrations in this paper required a lot of work. The  ``chaos game'', used in \cite{So4, BB} and for Figure \ref{bc}, is not sufficiently accurate. For the other figures, a Markov chain model was used, where $[0,1]$ was divided into $n$ intervals of equal length, with $4\cdot 10^3 \le n\le 4\cdot 10^6.$ For large magnifications, inverse iteration can be used, that is, we can directly count the offspring of points $x$ under the iteration of $G,$ as given in Definition \ref{bds} below. The drawback of that method is the increase of rounding errors under iteration of expanding maps. We do not discuss  numerical methods here since they serve for illustration only. 
 
Nevertheless, our opinion is that computer experiments and explicit approximation of fractal measures are appropriate tools for developing intuition and fostering communication. Bernoulli convolutions are a good testbed for fractal measures: on one hand they include a wide variety of examples, see \cite[Section 2]{ban}, and on the other hand they are ``pure type'' \cite[Section 3]{PSS} and have the same properties on each subinterval of $[0,1],$ due to their self-similar construction. Let us mention that some types of everyday ``big data'' come as fractal measures with several millions of bins: measurements of workload of a server, river discharges, precipitation, particulate matter, activity on financial markets. Singularities, that means extreme values, are of paramount importance in all these fields. \smallskip

The intuition given by Figure \ref{drei} is as follows. The upper function shows the typical BC, the appearance which by Solomyak holds for almost every $\beta .$ The function is probably continuous, but fractal. It shows more, but finer detail on higher resolution, like Weierstrass functions or Brownian motion. This could be made precise by determining a H\"older exponent or a dimension of the graph. For $\beta\to 2$ or $t\to\frac12$ we can expect that the H\"older exponent goes to zero and the dimension to 2. This increasing roughness is the most obvious trend which can be seen in BCs, see \cite[Section 2]{ban}. It was predicted from the convolution structure already by Erd\"os and others in the 1930s.

The middle panel of Figure \ref{drei} shows a singular measure, not a function. This is a typical BC of a Pisot number.
The variation of the approximating function is large, increased detail on higher resolution does not become finer, and values tend both to zero and infinity, yet not as fast as one might expect. The lower function represents a Perron parameter discussed in Section 10. This is a nowhere continuous function, with a dense set of poles, as we will show. The peaks of the function are higher than those of the comparable Pisot parameter, nevertheless it seems less singular than the Pisot BC. We do not know whether this BC is a singular measure, but we would guess that it is rather an $L_1$-function. Thus there are certain Perron parameters which lead to a lower degree of singularity than Pisot numbers.\smallskip

We have two conjectures. First, $\nu$ has a bounded density unless $\beta$ is weak Perron. This would imply that only countably many BCs are singular, and all rational parameters lead to bounded densities. Already Garsia \cite{G} made a similar conjecture, and most colleagues will probably agree.  The second conjecture, supported by our computer experiments, sounds more unlikely: the parameters $t$ for which $\nu_t$ is singular are nowhere dense. It could be even possible that $\Phi$ is continuous outside a nowhere dense set of parameters.\smallskip

We tried hard to make our presentation self-contained, lucid, and comprehensible to a wider audience.
The following section discusses the fractal construction of the parametric family of BCs. In Section 3 we introduce our type of multivalued map, in Section 4 we prove fundamental properties of finite orbits. In Section 5 we give examples of network-like orbits for the Fibonacci BC. In Theorem \ref{fiborbits} an infinite family of irreducible orbits of this type are constructed. In Section 6 we prove that
network-like orbits in the BC family are only possible for weak Perron parameters.  In Sections 7 and 8, we define address curves and discuss points with unique addresses. Section 9 and 10 deal with the intersections of address curves inside the overlap region. \bigskip

{\bf Acknowledgement.} The first version of this paper was written together with R\"udiger Zeller. The definition of branching dynamical system and the examples of Figure \ref{g75} and \ref{gm} are taken from his PhD dissertation \cite{Ze} which he defended in June 2015. Afterwards, he worked in other fields and finally drew back his authorship. I thank him for his cooperation, for many discussions and encouragement. Without R\"udiger, this paper would not exist.

I am very grateful to the referee for his careful review and a lot of helpful suggestions, in particular for improvements of Proposition \ref{grow} and Theorem \ref{dime}.

\section{Bernoulli convolutions as self-similar measures}\label{para}
We recall the fractal construction of the Bernoulli measure $\nu_t$ for fixed $t,$ as described in \cite{bar,BI,fal,hut} in a more general setting. We use the contracting similitudes $f_0,f_1$ which are the inverse maps of $g_0,g_1$ in \eqref{expa}.  
\begin{equation} 
f_0(x)=tx  \quad\mbox{ and } \quad  f_1(x)=tx+1-t  \quad \mbox{ on } I=[0,1] , \label{cont}\end{equation}

In the first step we consider the intervals $I_0=f_0(I)=[0,t]$ and $I_1=f_1(I)=[1-t,1]$ and the probability density $\phi^1(x)=\frac{1}{2t}(1_{I_0}(x)+1_{I_1}(x))$ which is piecewise constant. For $t>\frac12 ,$ there is a overlap interval $D=I_0\cap I_1=[1-t, t]$ on which $\phi^1$ assumes the larger value $\frac1t .$ In the second step, we consider intervals $I_{00}, I_{01}=f_0f_1(I), I_{10},I_{11},$ and we have new overlap intervals $D_0=f_0(D)=I_{00}\cap I_{01}$ and $D_1=f_1(D).$ 

In step $k,$ we consider all $2^k$ words $w=w_1...w_k\in \{ 0,1\}^k$ of length $|w|=k,$ contractions 
\begin{equation}\label{fwgw}
f_w=f_{w_1}\cdots f_{w_k}\, ,\quad\mbox{that is,}\quad f_w(x)=t^kx+(1-t)\cdot \sum_{j=1}^{k} w_j t^{j-1}\, ,
\end{equation}
and corresponding intervals $I_w=f_w(I)$ of length $t^k.$ The approximating density of $\nu_t$ on level $k$ can be written as 
\[ \phi^k(x)=(2t)^{-k}\sum_{|w|=k} 1_{I_w}(x) \ .\] 
That is, we count how many intervals $I_w$ contain the point $x.$  The normalizing factor $(2t)^{-k}$ is needed to obtain a probability measure. Hutchinson has shown that the measures given by the densities $\phi^k$ converge to the Bernoulli convolution $\nu_t$ in the space of probability measures with a metric which nowadays is called transport distance. See \cite{bar, fal, hut} for details.

Another construction uses the space $\Sigma=\{ 0,1\}^\infty$ of 01-sequences $s=s_1s_2...$ For each sequence $s$ there is a point $x$ in $[0,1]$ which can be written as 
\begin{equation}\label{pipi}
x=\pi(s)=\lim_{k\to\infty} f_{s_1}\cdots f_{s_k}(0)=\bigcap_{k=1}^\infty I_{s_1...s_k}\ .\end{equation} 
The mapping $\pi :\Sigma\to [0,1]$ is continuous \cite{bar,hut}. 
We say that $x$ has {\it address} $s=s_1s_2...$ In an overlapping system, a point can have many addresses. Let $\mu$ denote the product measure on $\Sigma$ given by $p_0=p_1=\frac12 .$ In probability theory, $\mu$ models a sequence of independent coin tosses. Here, it is considered the equidistribution of addresses. The Bernoulli convolution $\nu_t$ measures how many addresses the points $x\in [0,1]$ have: $\nu_t=\mu\cdot \pi^{-1}$ is the image measure of $\mu$ under the address map $\pi .$

With the concept of address, and the discussion of branching system which follows in the next section, the approximation $\phi^k(x)$ of $\nu_t$ has three interpretations: 
\begin{itemize} \item the number of intervals $I_w$ on level $k$ which contain $x,$
\item the number of prefixes of length $k$ of addresses of $x,$
\item the number of successors of $x$ in generation $k$ with respect to the system $G=\{ g_0,g_1\} ,$
 \end{itemize}
multiplied  by the normalizing factor $(2t)^{-k}.$ Convergence of
the density functions $\phi^k(x)$ to $\nu_t$ in the transport distance $d$ is fast, $d(\phi^k,\nu_t)\le t^k\cdot d(\phi^0,\nu_t)$ where $\phi^0$ is the equidistribution on $[0,1].$ However, $\nu_t$ need not be given by a density function.

These facts hold for all $t$ between 0 and 1, and we shall study $\nu_t$ for all overlapping cases $t\in [\frac12, 1)$ together. It will turn out that basic features of the complete family of Bernoulli convolutions are easier to understand than single specimen like Figure \ref{bc}.
As known from other parametric families of dynamical systems, e.g. the logistic family and the Mandelbrot set \cite{CE,MT,PR}, there is a lot of repetition and number theory which generates some structure in an otherwise chaotic picture. In our case, the number theory involves algebraic integers, mainly Pisot numbers, and address curves of periodic and preperiodic addresses will be the structural elements.  In Figure \ref{B5} the Bernoulli measures $\nu_t$ for $0.5\le t\le 0.76$ were represented by approximate densities $\phi^k$ with $k=26,$ which themselves were approximated by histograms with 4000 bins. 

\begin{figure}[h]
\includegraphics[width=0.95\textwidth]{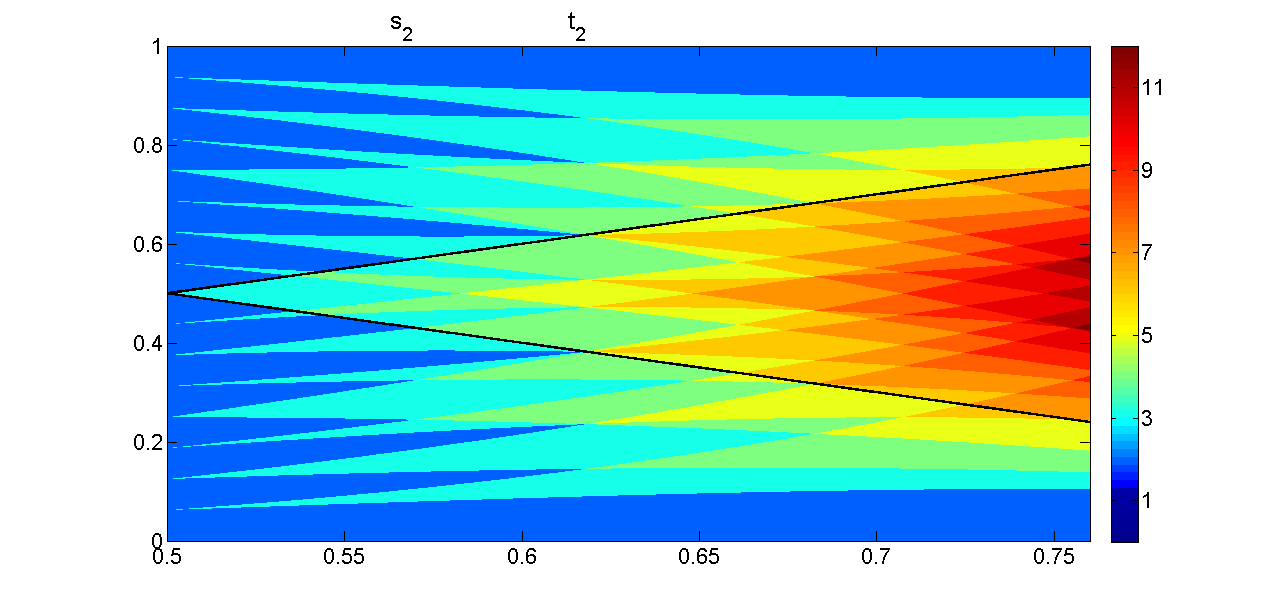}
\caption{Low-level approximation $\phi_t(x)$ of Figure \ref{B5}. For each parameter $t$ in $[0.5, 0.76],$ each point $x$ is assigned the number of pieces $I_w$ on level $k=4$ which contain $x,$ shown by color code. The black lines $x=t$ and $x=1-t$ mark the border of the main overlap region {\bf D}. }
\label{stufe4}
\end{figure}

To understand the measures better, we show in Figure \ref{stufe4} the function $\phi_t^k(x)$ on the smaller level $k=4.$ 
The normalizing factor is omitted to clarify the combinatorial structure. Possible values are the integers 1,...,16. For $t$ near to 1, all 16 intervals will contain the point $x=\frac12 ,$ but for $t\le 0.76$ at most 11 intervals intersect each other.  The number of intersections decreases when $t$ becomes smaller. For $t<\frac12$ all intervals $I_w$ within the same level are disjoint.
 For $x$ near to the endpoints of $[0,1],$ the value of $\phi^k$ is 1 since 0 is contained only in $I_{0...0}$ and 1 only in $I_{1...1}.$ On the other hand, we have the big triangular overlap region  ${\bf D}=\{ (t,x)|\ t\le x\le 1-t\}$ where two intervals meet already on the first level. We have other `horns' ${\bf D}_0=f_0({\bf D}), {\bf D}_w=f_w({\bf D})$ where intervals meet on levels 2,3,4. More rigorously, we  define 
\[ {\bf D}_w=\{ (t,x)|\ f_w(t)\le x\le f_w(1-t)\} \, ,\]
and note that this is a curvilinear triangle with tip at $(\frac12 , f_w(\frac12))$ and baseline $\{ (1,x)|\ 0\le x\le 1\}$ outside Figure \ref{stufe4}. Where horns intersect, the value of $\phi^k$ increases accordingly. This mainly happens for large $t$ and for points $x$ near $\frac12 .$
One could think that the maximum of $\phi^k(x)$ is always at $x=\frac12$ but for many $t$ this is not the case. The combinatorics of intersections of horns seems complicated.

While the intersections of horns indicate regions with many addresses, there are dark regions with baseline on the left border $t=\frac12$ in Figure \ref{stufe4}  which represent the points with only one address. It is interesting that their tips are at the golden mean parameter $t_2=\frac1\tau \approx 0.618$ with $\tau=\frac12(\sqrt{5}+1)\approx 1.618$ or smaller $t,$ like $t=s_2\approx 0.57 .$ When the level $k$ grows, these regions must become smaller, but they will not disappear completely.

Let us return to level 26, Figure \ref{B5}. Horns $f_w({\bf D})$ are still visible, but the dark regions on the left have become very thin. The golden mean $t_2\approx 0.618$ apparently marks a phase transition. For $t<t_2 ,$ there is a lot of structure. For $t>t_2 ,$ this structure becomes blurred, and only few details are apparent.
The interval for $t$ was chosen to include all inverses of Pisot numbers $\beta\in [1,2].$ The first and second Pisot numbers, 1.3247 and 1.3803, can be recognized  in Figure \ref{B5} at $t\approx 0.755$ and 0.725, respectively.  The next two Pisot numbers are vaguely visible at 0.693 and 0.682. 

\begin{figure}[h!]
\includegraphics[width=0.999\textwidth]{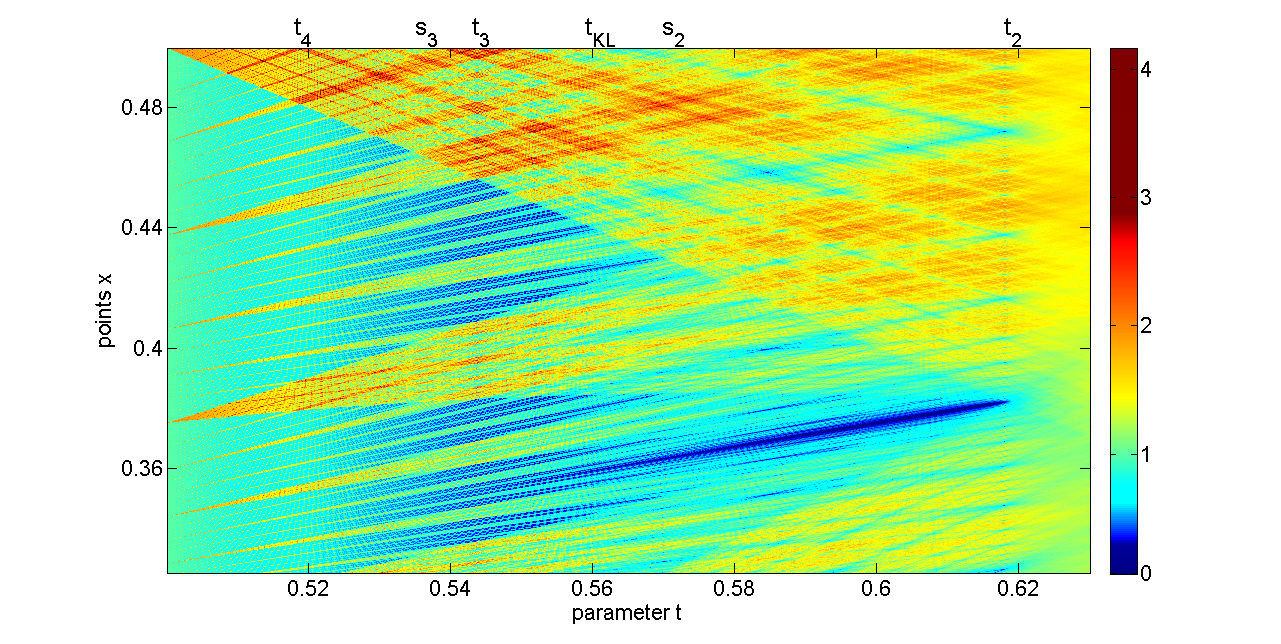}
\caption{$\Phi$ for $0.5\le t\le 0.63,\quad 0.325\le x\le 0.5.$ Right of the golden mean, the structure becomes blurred. On the left, periodic points appear outside {\bf D}, starting at $t_2\approx 0.618$ with the dark curve which hits the $y$-axis at $1/3.$}
\label{B6}
\end{figure}

In this note, we shall focus on the region $\frac12\le t\le t_2$ where the structure seems most apparent. Because of the symmetry of the measures $\nu_t$ we restrict our study to points $x\le \frac12 .$ Repetitions in Figure \ref{B5} at the lower and upper border are caused by the definition of $\nu_t.$  In \eqref{bern} one term on the right disappears when $A$ does not intersect the overlap interval $D$ since either $g_0$ or $g_1$ do not apply to $A.$  We have
\begin{equation} 
\nu_t(g_0(A))=2\nu_t (A)\ \mbox{ for }\ A\subseteq [0,1-t]\quad\mbox{ and }\quad
\nu_t(g_1(A))=2\nu_t (A)\ \mbox{ for }\ A\subseteq [t,1]\ .
\label{bern1} \end{equation}
For this reason, Figure \ref{B6} was restricted to $\frac12\le t\le 0.63$, $0.325\le x\le  0.5\, .$ 
Again, each $\nu_t$ was approximated by a histogram of the density $\phi_t^{26 }$  with 4000 bins, standardized so that the integral over $[0,1]$ is 1. As shown in the scale bar on the right, numerical values of densities on this approximation level vary only between zero and four, even though some $\nu_t$ are known to be singular.

In  Figure \ref{B6} we see again that for $t>t_2$ the structure becomes blurred. The same is true for $t$ near $\frac12$ where the $\nu_t$ converge to the uniform distribution. The most obvious feature are the dark curves which lead from the left border upwards and end just outside the overlap region {\bf D}, which here is the upper triangle with large density values. These curves indicate points with unique address which have been studied by many authors and will be dealt with in Section 7. In Sections 9 and 10, we shall consider the blue and red crossing points inside the overlap region which are related to singularity of the measures $\nu_t.$ Periodic points will play a key role in our investigation. 

Since pieces $I_w$ for fixed length $k=|w|$ can overlap, they can also coincide, which is called {\it exact overlap,} and has been considered a main reason for singularity of $\nu_t$ \cite{G}. By \eqref{fwgw}, $I_w=I_v$ holds if and only if
$f_w=f_v$ which is a polynomial equation for $t.$ Thus singularities - dark red spots in Figure \ref{B6} - are likely to appear at algebraic numbers $t$ and $\beta .$   Equation \eqref{fwgw} also says that maps $f_v,f_w$ with $|v|=|w|=k$ coincide iff   $f_w(0)=f_v(0)$ iff the fixed points coincide, that is $f_w(0)/(1-t^k)=f_v(0)/(1-t^k).$

\section{Branching systems}
Given an address $s=s_1s_2...$ we found the corresponding point $x=\pi(s)$ with contractions $f_{j}$ in \eqref{pipi}.
For fractal constructions where overlap of pieces can be neglected, at least in a measure-theoretic sense, for instance for $\nu_{\frac12},$ the mapping $\pi$ is essentially one-to-one. It is well-known that Bernoulli's product model of a sequence of independent coin tosses is isomorphic to the action of the doubling map $g(x)=2x\mbox{ mod }1$ on $[0,1].$ In the presence of substantial overlap, however, many points have a continuum of addresses. In such cases, one should better work with the expansive maps $g_j=f_j^{-1},$  as done in \cite{FS}.
We point out how this leads to a deterministic model of a branching process. 

Given a point $x$ in $[0,1]$ we are looking for {\it all} its addresses $w=w_1....w_k\in\{ 0,1\}^k$ for some level $k.$ To this end, we apply the maps $g_0(x)=f_0^{-1}(x)=\beta x$ and $g_1(x)=f_1^{-1}(x)=\beta x +1-\beta$ of \eqref{expa} repeatedly. In other words, we apply the mappings
\begin{equation}\label{gwfw}
g_w:I_w\to I\ \mbox{ with }\ g_w=g_{w_k}\cdots g_{w_1}\, ,\quad\mbox{that is,}\quad g_w(x)=\beta^kx+(1-\beta )\cdot \sum_{j=1}^{k} w_j \beta^{k-j}
\end{equation}
for all words $w$ with length $|w|=k.$ In contrast to the $f_w,$ the $g_w$ are only partially defined. Thus we can apply only those $g_w$ for which $x$ is in $I_w,$ and the corresponding words $w$ are just the addresses of $x,$ up to level $k.$

Note that the order of factors in $g_w=g_{w_k}\cdots g_{w_1}$ is the reverse of the order in $f_w=f_{w_1}\cdots f_{w_k}.$ The first letter of an address of $x$ is obtained by the first application of a map $g$ while it is determined by the last application of a map $f$ in a word $f_w.$ This convention will hold throughout the paper.

In this section, we formally present a concept of deterministic branching system which seems interesting in its own right.
It could be defined on arbitrary sets - for instance on finite sets, with a combinatorial flavor as in Proposition \ref{grow} below. For the purpose of this paper, however, we define it on intervals. We also require that each point $x$ is in the domain of some $g_i$ so that `survival of the branching process is guaranteed'.  This excludes Bernoulli convolutions with $\beta >2,$ or $t<\frac12 ,$ where the measure lives on a Cantor set. Roughly speaking, we consider multivalued maps on real intervals which consist of a finite number of branch maps, as in \cite{I5}. 

\begin{Definition}[Branching dynamical system on an interval]\label{bds}\hfill\\
A branching dynamical system on $[a,b]\subset\RR$ is given by a set of subintervals $I_i\subset [a,b]$ with $\bigcup_{i=1}^m I_i=[a,b],$ and continuous mappings $g_i:I_i \rightarrow [a,b], $  $i=1,\hdots, m.$\\
For each point $x\in [a,b], $ the set of immediate successors is $G(x)=G^1(x)=\{g_i(x)\ |\ x\in I_i\}$. The successors of $q^{th}$ generation are given by the recursion $G^q(x)=G(G^{q-1}(x))$.
The orbit of $x$ is $O(x)=\bigcup_{q=0}^\infty G^q (x)$ with $G^0(x)=\{ x\} .$
\end{Definition}

The orbits of a branching dynamical system have the structure of a rooted infinite tree, with 1 up to $m$ immediate successors for each vertex. If some values in the orbit coincide, the tree structure is modified and cycles appear. Here we are interested in the case that various successors coincide and the tree becomes a finite network. 

An interesting problem for branching dynamical systems is the existence of an absolutely continuous invariant measure $\mu$ with a corresponding growth factor $\lambda$ in the sense that

\begin{equation} \lambda\cdot \mu (A)=\sum_{i=1}^m \mu(g_i^{-1}(A))\quad \mbox{ for Borel sets } A\subset [a,b] .
\label{invme}\end{equation}

Even in the non-overlapping case where the $I_i$ touch at endpoints, the existence of $\mu$ could only be proved when the $g_i$ are assumed to be monotone, differentiable and expansive, see \cite{LM}.
Next, one could ask for ergodic theorems, which involves the question of asymptotic independence of successors on disjoint branches of the orbit of $x.$  One may also ask for the growth factor $\lambda_x$ of the number of successors for different $x,$ and its connection with $\lambda ,$  see \cite{FS,Kempton} and Section 4. It would be very interesting to know under which conditions the $G^q(x)$ approach an asymptotic distribution for $q\to \infty .$
\ These problems are difficult, so we consider more special systems. 

\begin{Definition}[Expansive and linear branching systems]\label{fbs} \hfill\\
A branching dynamical system $\{ g_i:I_i \rightarrow [a,b], $  $i=1,\hdots, m\}$  is called expansive if
each $g_i$ maps onto $[a,b],$ and is monotonuous with derivative $|g'_i(x)|>\gamma>1$ for all $x\in I_i.$\\
An expansive branching system is called linear if $g_i(x)=\beta_i x+z_i$ are linear functions. If all $\beta_i$ are equal to $\beta ,$ we have a linear branching system with slope $\beta .$ A linear branching system with $m=2$ linear functions with slope $\beta $ will be called Bernoulli system or Bernoulli convolution.
\end{Definition}

For an expansive branching system, the inverse functions $f_i=g_i^{-1}:[a,b]\to I_i$ are contractive. The family $\{ f_1,...,f_m\}$ then forms an iterated function system with attractor $[a,b],$ cf. \cite{bar,fal}.  For any choice of probabilities $p_1,...,p_m\ge 0$ with $\sum p_i=1$ there is a unique probability measure $\nu$ which fulfils the equation 
\begin{equation} 
\nu (A)=\sum_{i=1}^m p_i\nu (f_i^{-1}(A))= \sum_{i=1}^m p_i\nu(g_i(A))
\quad\mbox{ for Borel sets } A\subset [a,b] . \label{self}\end{equation}
For a linear branching system,  Lebesgue measure can be taken as the invariant measure $\mu$ in \eqref{invme}, with growth factor $\sum 1/\beta_i .$  If all maps have slope $\beta ,$ this is $\lambda=\frac{m}{\beta}.$ 

For a Bernoulli system, $\beta$ must be between 1 and 2, and the endpoints of the interval are the fixed points of the two maps which are denoted $g_0$ and $g_1$ in this case. See Figure \ref{bc}. For each $\beta$ there is only one conjugacy class of these branching dynamical systems. To obtain the same basic interval for different values of $\beta ,$ we take $[0,1]$ and the maps $g_0,g_1$ in \eqref{expa}. 

In the 1930s, the term \emph{Bernoulli convolution} was used for the self-similar measure $\nu$ corresponding to $\{ f_0,f_1\}$ in \eqref{cont} and $p_1=p_2=\frac12 ,$  the special case \eqref{bern} of equation \eqref{self}, including the non-overlapping case $\beta >2$. 
The name comes from the fact that $\nu $ (in the version $g_0(x)=\beta x, \ g_1(x)=\beta x-1$ on $[0,\frac{1}{\beta-1}]$), can be expressed as an infinite convolution
$\nu =\bigotimes_{n=1}^\infty \nu_n$ of two-point measures $\nu_n=\frac12\delta_0+\frac12\delta_{t^n}$ where $t=1/\beta .$ In probabilistic terms, $\nu$ is the distribution of a random sum $\sum_n \omega_n t^n$ where all coefficients $\omega_n$ are independently chosen as $0$ or $1$ with probability $\frac12 .$ 
There are two reasons to focus on this special case. On the one hand, even Bernoulli convolutions are not well understood. On the other hand,  we hope that results for the simple system will provide insight into more general cases.

\section{Basic properties of finite orbits}
The measure $\nu$ will now be studied for different $\beta$ with the help of finite orbits. First we show that each finite orbit determines the measure $\nu (B)$ of certain intervals $B.$  

\begin{Proposition}[Markov partition generated by a finite orbit] \label{mark}\hfill\\
Let $\{ g_i:I_i\to [a,b]\, |\, i=1,...,m\}$ be an expansive branching system, and let $X$ be a finite set which contains its successors: $g_i(X)\subseteq X$ for all $i$. Let $(p_1,...,p_m)$ be a vector of positive probabilities and $\nu$ the corresponding self-similar measure defined by (\ref{self}). Then the points of $X$ define a partition of $[a,b]$ into subintervals $J_k$ for which $\nu (J_k)$ is uniquely determined. 
\end{Proposition}

{\it Proof. }  Definition \ref{fbs} implies that the successors of $a$ and $b$ are either $a$ or $b.$ So we add these two points to $X$ without changing the assumptions. The points of $X$ will be ordered, $a=x_0<x_1<...<x_n=b. $ The partition intervals are $J_k=[ x_{k-1},x_k]$ for $k=1,...,n.$  For each $i$ and $k,$ the set $g_i(J_k)$ is either empty or a union of consecutive intervals $J_\ell .$ (Since $g_i$ is monotonic on $I_i,$ so in the case $x_{k-1}\in I_i, x_k\not\in I_i$ either
$g_i(J_k)=[g_i(x_{k-1}),b]$ or $g_i(J_k)=[a,g_i(x_{k-1})].$ The other cases are similar.)\\
We now define an $n\times n$-matrix 
\[ M=(m_{k\ell})\quad \mbox{ with } \quad m_{k\ell}=\sum \{p_i\, | \, J_\ell\subseteq g_i(J_k)\} \]   
where $m_{k\ell}=0$ when no $i$ fulfils the condition.  The matrix was constructed in such a way that the condition (\ref{self}) for the measures $w_k=\nu (J_k)$ of the sets $A=J_k, k=1,...,n$ coincides with the matrix equation $w=Mw,$ that is, $w_k=\sum_\ell m_{k\ell}w_\ell$ for $k=1,...,n.$

Each column sum of $M$ equals $\sum p_i =1$ because $J_\ell$ is once  in 
$g_i(I_i)$ for each $i.$ In other words, the vector $(1,...,1)$ is a left eigenvector of $M$ for the Perron-Frobenius eigenvalue 1. So there is a non-negative right eigenvector $w=(w_1,...,w_n)'$ of $M,$ that is, a solution of $Mw=w.$ 
To verify the uniqueness of $w,$ we have to show that the matrix $M$ is irreducible.

Given $k,\ell\in\{ 1,...,n\},$ we have to find a finite sequence $k_1,...,k_q\in\{ 1,...,n\}$ with $m_{k_{j-1}k_j}>0$ for $j=1,...,q,$ with $k_0=k, k_q=\ell .$ In terms of partition intervals, this is 
\[ J_{k_1}\subseteq g_{i_1}(J_k),\,  J_{k_2}\subseteq g_{i_2}(J_{k_1})\ ,..., \  J_\ell=J_{k_q}\subseteq g_{i_q}(J_{k_{q-1}})\ .\] 
Now since the inverse maps $f_i$ of the $g_i$ are contractive, there is a word $i_1i_2...i_q$ with $f_{i_1}f_{i_2}...f_{i_q}(J_\ell)\subseteq J_k$ (cf. \cite{bar,fal}). In other words,  $J_\ell\subseteq g_{i_q}...g_{i_2}g_{i_1}(J_k)\, .$ Since the $g_i$ are strictly monotonic, $J=g_{i_{q-1}}...g_{i_1}(J_k)$ is an interval with endpoints in $X,$ and we can find a subinterval
$J_{k_{q-1}}$ with $J_\ell\subseteq g_{i_q}(J_{k_{q-1}})\subseteq g_{i_q}(J).$ Recursively, the other $J_{k_j}$ are defined. \hfill $\Box $\smallskip

{\it Remarks.} The condition $g_i(X)\subseteq X$ is quite strong. Theorem \ref{nonex} below implies that for Bernoulli convolutions with $t>t_2,$ sets $X\not\subseteq \{ a,b\}$ exist only if $\beta$ is a weak Perron number. For arbitrary $t\le t_2$ such sets are given in Section 7, but they will not intersect the overlap region $D.$

If the $p_i$ are rational, then the coefficients of the equation $Mw=w$ and the measure values $w_k=\nu (J_k)$ will also be rational. Moreover,  in the case of Bernoulli convolutions the constructed irreducible Markov chain is non-periodic since $g_0(0)=0$ implies $m_{11}=p_1=\frac12 .$  \medskip

Next, we study the growth of successor generations for a point $x$ with finite orbit and the local dimension of $\nu$ at $x.$  Proposition \ref{grow} below generalizes Theorem 5.2 of Baker \cite{Baker}. Under special assumptions,  Feng \cite{F4}, Feng and Sidorov \cite{FS}, 
and  Kempton \cite{Kempton} have found equations or inequalities related to Propositions \ref{grow} and Theorem \ref{dime} which hold for all $x,$ or for Lebesgue almost all $x.$ Here we consider finite sets
$X=\{ x_1,...,x_n \}$ in a general setting.

Let us consider the $n\times n$ successor matrix 
\[ S=(s_{k\ell})\quad\mbox{ with }\quad s_{k\ell}= |\{ i\, |\, g_i(x_k)=x_\ell \} | \quad\mbox{ for } k,\ell =1,2,...,n\ .\]
$S$ is the adjacency matrix of a directed graph, with an edge from $x_k$ to $x_\ell$ for each $i$ with $g_i(x_k)=x_\ell .$ The sum of the $k$-th row is the number of successors of $x_k$ within $X.$ The matrix $S$ is irreducible if for each $x_j$ and $x_k,$ there is a path in the graph from $x_j$ to $x_k.$

It is well-known and easy to check that for such an adjacency matrix, the entry $s^q_{k\ell}$ of the matrix power $S^q$ counts the number of paths of lengths $q$ from 
$x_k$ to $x_\ell ,$ for $q=1,2,3,...$ Thus $s^q_{k\ell}$ is the number of compositions $g=g_{i_1}...g_{i_q}$ of the $g_i$ which fulfil $g(x_k)=x_\ell .$ The $k$-th row sum of $S^q$ is the number $\gamma^q(x_k,X)$ of $q^{\rm th}$ generation successors of $x_k$ within $X.$ We write $\gamma^q(x,X)$ instead of $|G^q(x)\cap X|$ since points are counted with multiplicity. Different $g$ may lead to the same point. Let $\gamma^q(x)$ denote the number of all $q^{\rm th}$ generation successors of $x,$ counted with multiplicity. If $\lambda_x=\lim_{q\to\infty} \sqrt[q]{\gamma^q(x)}$ exists, it is called the growth factor of $x.$ If the limit does not exist, liminf and limsup will be called lower growth factor
$\underline{\lambda}_x$ and upper growth factor $\overline{\lambda}_x$ of $x.$

\begin{Proposition} [Growth factor of a finite set and finite orbit] \label{grow}\hfill\\
Let  $X=\{ x_1,...,x_n \}$ have the successor matrix $S$ in a branching dynamical system, and let $\rho$ denote the spectral radius of $S.$   
\begin{enumerate} 
\item[ (i)] The numbers $\gamma^q(x_k,X)$ of $q^{\rm th}$ generation successors of $x_1,...,x_n$ within $X$ are given by the vector  $S^q\cdot (1,...,1)' $ for $q=1,2,...$ 
\item[(ii)] There is a point $x$ in $X$ with \  $\underline{\lambda}_x \ge \rho\, .$ If $S$ is irreducible, this holds for all $x$ in $X.$ 
\item[(iii)]  If $X$ is the finite orbit of $x_1$ in the branching dynamical system, the growth factor $\lambda_{x_1}$ of $x_1$ exists and  equals $\rho .$ If $S$ is irreducible, then all $x_j$ have growth factor $\rho .$
\end{enumerate}
\end{Proposition} 

{\it Proof. }  (i) was proved above. \ (ii): Since $S$ is a non-negative matrix, there is a non-negative eigenvector $v=(v_1,...,v_n)'$ of $S$ with eigenvalue $\rho ,$ and $v$ is positive if $S$ is irreducible \cite{BP}. We can assume $v_k\le 1$ for all $k.$ For each $k$ with $v_k>0$ we have
\[ \gamma^q(x_k)\ge  \gamma^q(x_k,X) = (S^q\cdot (1,...,1)')_k \ge (S^q\cdot v)_k =\rho^q v_k \ .\]
Taking $q^{\rm th}$ root and liminf, we obtain (ii).

(iii): If $X$ is a finite orbit, then $g_i(X)\subseteq X$ for $i=1,...,m,$ so  $ \gamma^q(x_k,X)= \gamma^q(x_k)$ for all $k.$ By (i), the row sum norm of $S^q$ is $\| S^q\| =\max_{j=1}^n \gamma^q(x_j)\, .$  Since each $x_j\in O(x)$ has the form $x_j=g_{i_1}...g_{i_r}(x_1)$ with $r\le n ,$ we have
$ \gamma^{q-n}(x_j)\le \gamma^q(x_1)$ for $q>n.$ This implies 
\begin{equation}
 \gamma^q(x)\le \| S^q\| \le m^n\cdot \gamma^q(x)\quad \mbox{ for } q>n\, . \label{norm}
\end{equation}
Since the spectral radius of $S$ fulfils $\rho =\lim_{q\to\infty}\sqrt[q]{\| S^q\| }$ for every matrix norm, this proves  $\lambda_{x_1}=\rho .$ If $S$ is irreducible, the same argument works when $x_1$ is replaced by any $x_j.$
 \hfill $\Box \vspace{2ex}$

For an expansive branching system and a self-similar measure $\nu ,$ the growth factor can be expressed in terms of local dimension. The {\it local dimension} of a measure $\nu$ on $\mathbb R$ at a point $x$ is defined as 
\begin{equation} \label{dimdef}
d_x(\nu )=\lim_{s\to 0}\frac{\log \nu(U(x,s))}{\log s}\qquad \mbox{ where }\quad U(x,s)=[x-s,x+s]\ .
\end{equation}
If the limit does not exist, liminf and limsup will be called lower local dimension 
$\underline{d}_x(\nu)$ and upper local dimension $\overline{d}_x(\nu)$ of $\nu$ at $x,$ respectively. The following theorem will be used to detect poles and zeros in Bernoulli measures, as explained at the end of Section 5. The referee noted that it can be stated more clearly for general self-similar measures \eqref{self} with probabilities $p_1,...,p_m.$ 
This could be used to investigate biased Bernoulli convolutions as studied in \cite{JSS}. To this end, we define
a weighted successor matrix $\tilde{S}$ for the set $\{ x_1,...,x_n \}$ as follows. 
\[ \tilde{S}=(\tilde{s}_{k\ell})\quad\mbox{ with }\quad \tilde{s}_{k\ell}= \sum\{ p_i\, |\, g_i(x_k)=x_\ell \}  \quad\mbox{ for } k,\ell =1,2,...,n\ .\]
Actually, we should write $\tilde{S}_p$ but we suppress the index. For the case of equal probabilities $p_i=\frac1m $ we have  $\tilde{S}=\frac1m\cdot S .$ 

\begin{Theorem} [Local dimension at finite orbits in branching systems of slope $\beta$] \label{dime}\hfill\\
Let $\{ g_i:I_i\to [a,b]\, |\, i=1,...,m\}$ be a linear branching system with slope $\beta ,$ and 
$X=\{ x_1,...,x_n \}$ a finite subset of $[a,b].$ Let $p_1,...,p_m$ be positive with $\sum p_i =1,$ and $\tilde{S}$ the corresponding weighted successor matrix with spectral radius $\tilde{\rho}.$ Let $\nu$ be the self-similar measure generated by the $g_i$ and the  probabilities $p_i$ defined in \eqref{self}. Then 
\begin{enumerate} 
\item[ (i)]  $\overline{d}_x(\nu )\le \frac{-\log \tilde{\rho}}{\log \beta}$ for at least one $x$ in $X.$ If $S$ is irreducible, this holds for all $x$ in $X.$ 
\item[(ii)]  If $X$ is the orbit of $x_1$ in the branching dynamical system then 
$d_{x_1}(\nu )= \frac{-\log \tilde{\rho}}{\log \beta}\ .$ In the case of irreducible $S$ all 
$x_j$ have this local dimension. For equal probabilities $p_i=\frac1m$ we have $d_{x_1}(\nu )= \frac{\log m -\log \rho}{\log \beta}\ .$
\end{enumerate}
\end{Theorem} 

{\it Proof. }  First we note that equation \eqref{self}, iterated $q$ times, yields the same equation for words $w=w_1...w_q$ of length $q.$
\begin{equation}\label{selfw} 
\nu (A)=\sum_{|w|=q} p_w\nu (f_w^{-1}(A))= \sum_{|w|=q} p_w\nu(g_w(A))
\quad\mbox{ for Borel sets } A\subset [a,b] . 
\end{equation}
We apply this to the set $A=U(x_k,\alpha\beta^{-q})$ where $\alpha$ is a positive number. By \eqref{gwfw} each $g_w$ is a similitude with factor $\beta^q.$ Thus $g_w(A)=U(g_w(x_k), \alpha).$ Terms with $g_w(x_k)$ outside $X$ are disregarded. The sum of $p_w$ with $|w|=q$ and  $g_w(x_k)=x_\ell$ is just the term $\tilde{s}_{k\ell}^q$ of the matrix power $\tilde{S}^q$ as can be verified by induction on $q.$ Let $c=\min_{\ell=1}^n \nu(U(x_\ell , \alpha)).$ Then
\begin{equation}\label{nuab} 
\nu(U(x_k,\alpha\beta^{-q})) \ge \sum_{\ell=1}^n \tilde{s}_{k\ell}^q \nu(U(x_\ell , \alpha))
\ge c\cdot \sum_{\ell=1}^n \tilde{s}_{k\ell}^q = c\cdot (\tilde{S}^q\cdot(1,...,1)')_k\ .
\end{equation}
The non-negative matrix $\tilde{S}$ has a non-negative eigenvector $v=(v_1,...,v_n)'$ with $v_k\le 1$ where the eigenvalue is the spectral radius $\tilde{\rho}.$ We choose an $x=x_k$ with $v_k>0$ so that 
\begin{equation}\label{nucd}
\nu(U(x,\alpha\beta^{-q})) \ge c\cdot (\tilde{S}^q\cdot(1,...,1)')_k\ge c\cdot (\tilde{S}^q\cdot v)_k =c\cdot \tilde{\rho}^q v_k\ . 
\end{equation}
To estimate local dimension at $x,$ we take logarithms and divide by $\log \alpha\beta^{-q}$ where $q$ is so large that the denominator is negative:
\[ \frac{\log \nu(U(x,\alpha\beta^{-q}))}{\log \alpha -q\log \beta}\le  \frac{\log cv_k +q\log\tilde{\rho}}{\log \alpha -q\log \beta} \]
When $q$ tends to infinite, the limit on the right-hand side is $-\log \tilde{\rho}/\log \beta ,$ independently of $\alpha .$
The upper limit of the left-hand side cannot be larger, for any $\alpha .$ This proves (i). \smallskip

For (ii) we assume that $X$ is the orbit of $x_1.$ On the one hand this means $g_i(X)\subseteq X$ so that there are no $g_w(x_k)$ outside $X$ and the first inequality of \eqref{nuab} is an equation. On the other hand, each $x_k$ has the form  $x_k=g_{w^k}(x_1)$ for some word $w^k$ of length $r_k<n.$ Equation \eqref{selfw} with $A=U(x_1,\alpha\beta^{-r_k})$ provides inequalities 
\[ \nu(U(x_1,\alpha))\ge  \nu(U(x_1,\alpha\beta^{-r_k})) \ge p_{w^k} \nu(U(x_k,\alpha)) \ .\]
Taking $C$ as the maximum of the values $ \nu(U(x_1,\alpha))/ p_{w^k}$ over $k=2,...,n$ and rewriting \eqref{nuab} we get
\[ \nu(U(x_1,\alpha\beta^{-q})) = \sum_{k=1}^n \tilde{s}_{1k}^q \nu(U(x_k , \alpha)) \le C \| \tilde{S}^q\|
\le 2C \tilde{\rho}^q \]
for sufficient large $q,$ where $ \| \tilde{S}^q\|$ denotes the row sum norm. Taking logs and dividing by $\log \alpha\beta^{-q}$ as above for \eqref{nucd} we conclude $\underline{d}_{x_1}(\nu)=-\log \tilde{\rho}/\log \beta ,$ If $S$ is irreducible, every $x_j$ can replace $x_1$ in this calculation. For probabilities $p_i=\frac1m $ we have $\tilde{\rho}=\frac{\rho}{m}.$    \hfill $\Box$ \medskip

{\it Remark.}  The one-dimensional structure of the space $[a,b]$ was used only in Proposition \ref{mark}. Proposition \ref{grow} is just a combinatorial statement. Theorem \ref{dime} holds in $\mathbb R^n $ for similarity maps with equal contraction factors, and balls $U(x,s)$ instead of intervals.

\section{Network orbits in the Pisot case}
\begin{Definition}[Different types of algebraic integers]\hfill\\
An algebraic number $\beta$ is a root of a polynomial
$p(x)=b_mx^m+b_{m-1}x^{m-1}+\hdots+b_1 x+b_0$ with integer coefficients. The degree $m$ of $p$ and $b_m>0$  are chosen to be minimal, so $p$ is unique and we can say $\beta$ has degree $m.$  If $b_m=1$ then $\beta$ is called algebraic integer. Here we consider real algebraic integers $\beta$ which are greater than one.
If all other roots (``conjugates'') $\lambda$ of $p$ have modulus strictly smaller than one, $\beta$ is called a Pisot number.
If only $|\lambda|\le 1$ is required and equality holds for one root, then equality will hold for all conjugates, except for one real root which is $1/\beta .$ In this case $\beta$ is called a Salem number.  If the conjugate roots only fulfil $|\lambda|<\beta$ we call $\beta$ a Perron number. If only $|\lambda|\le\beta$ is required, $\beta$ is a weak Perron number.  Finally, if all conjugates fulfil $|\lambda|>1$ and  $|b_0|=2,$ we call $\beta$  a Garsia number. 
\end{Definition}

{\bf Important Pisot parameters. } The best known Pisot number is the Fibonacci or golden number  $\tau=\frac12(\sqrt{5}+1)\approx 1.618\, .$ More generally, for each $n\ge 2$ there is a multinacci parameter $\tau_n$ defined as positive solution of $x^n=x^{n-1}+...+x+1.$ The parameters $t_n=1/\tau_n$ are the most obvious landmarks for the function $\Phi ,$ shown in many of our figures. 
There is another family $\varphi_n$ of Pisot numbers, defined as solutions $x>1$ of $x^{n+1}=2x^n-x+1$ for $n\ge 2.$ We call them doubling numbers and write $s_n=1/\varphi_n.$ These parameters are also indicated in our figures. Figure \ref{bc} refers to $\varphi_2.$ \medskip

The set of Pisot numbers is closed. The $\tau_n, \varphi_n$ for $n=2,3,...$ and a single number $\chi$ mentioned in Section \ref{geom} are the only accumulation points of Pisot numbers in $(1,2).$ See \cite{AF} for more details and references. There are still many mysteries concerning Pisot, Salem and Garsia numbers, cf. \cite{HP}. We hope that the study of Bernoulli convolutions can shed some light on their structure. 

For Garsia numbers $\beta\in (1,2)$ it is known that the Bernoulli convolution possesses a bounded density \cite{G}. For Pisot numbers, the Bernoulli convolutions are singular \cite{E}. For all other $\beta ,$ including all rational numbers and Salem numbers, this question is not yet resolved.  \medskip

The following `folklore theorem' says that for Pisot numbers $\beta$ there are lots of network-type orbits. Garsia \cite{G} found the basic property, a uniform discreteness lemma, which was recently used by Baker \cite{Baker} to verify the statement. Schmidt \cite{Sch} proved a special case with other estimates. A direct proof is given in the first arxiv version of this paper. 
		
\begin{Theorem}[Abundance of finite orbits for Pisot slope \cite{G,Sch}] \label{finiteorbits}\hfill\\
	Let $g_i(x)=\beta x+z_i$, $i=1,...,k$ be a branching dynamical system on an interval $[a,b].$ If
	$\beta$ is a Pisot number and each $z_i$ is in $\QQ(\beta),$
	then $O(x)$ is a finite set for each $x$ in $\QQ(\beta).$
\end{Theorem}

\begin{figure}[h]
\begin{center}
\includegraphics[width=0.6\textwidth]{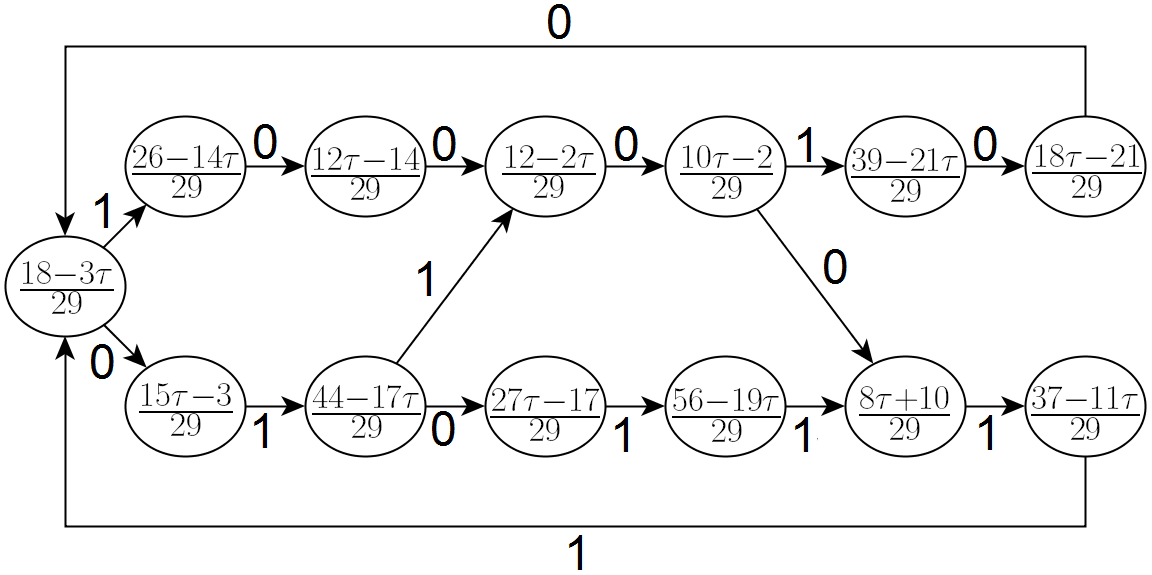}
\end{center}
\caption{The orbit of $\frac{18-3\tau}{29}$ for the golden mean case}
\label{g75}
\end{figure}

The Bernoulli convolution for the Fibonacci parameter was studied by many authors \cite{AZ,F4,FS,Hu,LP,SV}, with  focus on the dimension and multifractal spectrum of the measure. Here we consider periodic orbits for $g_0(x)=\tau x\, ,\ g_1(x)=\tau x+1-\tau$ to illustrate our subject.  
The point $\frac12$ has period 3, it fulfils $g_{100}(\frac12 )=g_{011}(\frac12 )=\frac12 .$ Thus the orbit $O(\frac12)$ consists of five points, and the growth factor is $\sqrt[3]{2}\approx 1.2599.$ By Theorem \ref{dime}, the local dimension of the Bernoulli measure at $\frac12$ is 0.9603. A slightly more complicated case  is shown in Figure \ref{g75}. The point $x_0=\frac{18-3\tau}{29}$ has the fixed period 7, and it has 5 successors in generation 7. The growth rate is $\sqrt[7]{5}\approx 1.2585 ,$ and the local dimension is 0.9626.

The networks of $\frac12$ and $x_0$  have a special structure. Let us say that a directed graph is a \emph{mixture of cycles} of period $p$ if all minimal directed cycles in the graph have length $p,$ and each edge is contained in a cycle. A vertex $x_0$ which is contained in all cycles of the mixture will be called \emph{root.} 
Let $m$ denote the number of cycles. If $x_0$ is a root, the number of successors of $x_0$ in generation $p\cdot k$ is $m^k$ for $k=1,2,....$ So the growth factor is $\lambda_{x_0} =\sqrt[p]{m}.$ 

Vertices in a network like Figure \ref{g75} are points in $[0,1],$ and edges are labelled with 0 or 1, depending on which action $g_i$ they represent.  Thus a cycle of $p$ edges starting at vertex $x$ is labelled by a 01-word $w=w_1...w_p$ and indicates that $g_w(x)=x,$ or $f_w(x)=x$ when we run through the cycle in reverse direction.  If another cycle of the same length starts at $x,$ with labelling $v=v_1...v_p,$ then $x$ is the common fixed point of $f_v$ and $f_w.$ According to the last paragraph of Section 2, this implies $f_v=f_w$ and exact overlap of the pieces: $I_v=I_w.$  So a mixture of $m$ cycles of period $p$ with root $x_0$ is just an illustration of the fact that $x_0$ is the common fixed point of mappings $f_{w^1},...,f_{w^m},$ which means that these mappings coincide, and there is exact overlap of the corresponding intervals: $f_{w^1}=...=f_{w^m}$ and $I_{w^1}=...=I_{w^m}.$

There is a simple product construction for cycle mixtures with root, based on the fact that  $f_{w^1}=...=f_{w^m}$ together with  $f_{v^1}=...=f_{v^n}$ implies the equality of all $f_{w^iv^j}$ for $i=1,...,m$ and $j=1,...,n.$ Graphically, the two graphs are concatenated by redirecting the incoming edges of one root to the other root so that cycles from a root to itself are now labelled $u=w^{i}v^{j}$ or $u=v^{j}w^{i},$ respectively. When the product graph with two roots is constructed, new values $x$ have to be computed for each vertex, from equations of the form $f_u(x)=x$ which for contractive maps have a unique solution. The growth of the product graph is between the growth of the two factors, as can be easily checked.

Another construction for the Fibonacci case will give an impression of how complicated finite orbits can be, even for cycles of equal length. Pisot parameters $\beta$ of higher degree, with conjugates nearer to the unit circle, admit much more intricate networks, cf. Section \ref{super}.
Let us say that a cycle mixture is \emph{prime} if it contains exactly one root vertex, and thus is not obtained by concatenation. 

\begin{figure}[h]
\begin{center} \
\includegraphics[width=0.46\textwidth]{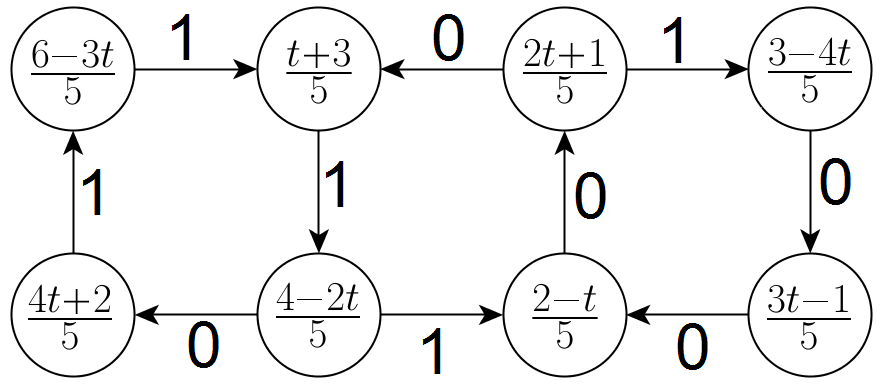}\hfill
\includegraphics[width=0.49\textwidth]{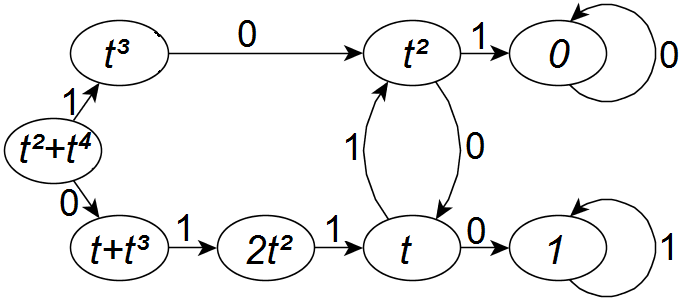}
\end{center}
\caption{Orbits with maximum and minimum growth rate}
\label{gm}
\end{figure}

\begin{Theorem}[Finite orbits with large complexity in the Fibonacci case] \label{fiborbits}\hfill\\
For each odd period $p=2k-1, k\ge 1,$ the golden Bernoulli convolution admits a prime cycle mixture of period $p$ with $n=4k-3$ vertices and $m=F_{k+1}$ cycles, where $F_1=F_2=1, F_3=2,...$ denote the Fibonacci numbers.  For $k\to\infty$ the growth rates $\rho_k=\sqrt[p]{m}$ of these orbits converge to the maximum growth rate $\sqrt{\tau} .$
This maximum growth rate is realized by the network on the left of Figure \ref{gm}, a mixture of 3 cycles of lperiod 4 without root.  	
\end{Theorem}

{\it Proof.} The cases $p=3$ and $p=7$ were discussed above, and we noted the basic equation $f_{100}=f_{011}$ which follows from the defining relation $1=t+t^2$ of the Fibonacci number $t=t_2.$ From the basic equation we  conclude $f_{u100}=f_{u011}$ and $f_{100v}=f_{011v}$  for all words $u,v.$  With $v=00$ and $u=01$ we get
$f_{10000}=f_{01100}=f_{01011}$ which settles the case $p=5.$  Repeating the step with $v=11$ and $u= 1000, 0110$ we get the case $p=7$ in Figure \ref{g75}.

Now we use induction to extend the graph further. Assume that for $k,$  we have $m=F_{k+1}$ words of length $p=2k-1$ describing the same map $f,$ where $F_k$ of the words have the form  $u0$ and $F_{k-1}$ have the form $v1,$ (this holds true for $k=3,$ in each further step we interchange 0 and 1 in this last statement). For the words ending with 0 we define two extensions $u011$ and $u100,$ for the words ending with 1 we define the extension $v111.$ All extensions represent the map $\overline{f}=f\cdot f_{11}.$ Their number is $2F_k+F_{k-1}=F_{k+2}.$ Moreover, $F_k+F_{k-1}=F_{k+1}$ of them end with 1, while $F_k$ end with 0. The new root is the fixed point $y$ of $\overline{f}.$ The extended graph has 4 new vertices which represent certain paths from the root. The former cycles $u0,v1,$ lead to the vertex $f_{11}(y),$ their extension $u01,v11$ lead to vertex $f_1(y),$ the paths $u1$ lead to $f_{00}(y)$ and  their extensions $u10$ to $f_0(y).$ Obviously, these vertices are not roots of the graph. The statement is shown for $k+1$ which completes the induction. 

The convergence $\sqrt[2k-1]{F_{k+1}}\to\sqrt{\tau}$ follows from $F_k\approx \tau^k/\sqrt{5}.$ The minimum local dimension of $\nu_t$ for the Fibonacci case was found by Hu \cite{Hu}, cf. \cite{F4}. The corresponding growth rate $\sqrt{\tau}\approx 1.2720$ is realized by the orbit of $x=\frac{2\tau -1}{5}=1/\sqrt{5}$ shown on the left of Figure \ref{gm}. One image of $x$ goes to $1-x$ after two steps, the other returns to $x$ after four generations. The growth rate $\rho$ fulfils the equation $\rho^4=\rho^2 +1$ (due to symmetry, the spectral radius can be determined for a graph with 4 vertices). Thus $\rho =\sqrt{\tau},$ with local dimension 0.9404. In Section \ref{super} we shall see that such double orbits combining $x$ with $1-x$ often yield large growth rates.  We briefly prove that the growth rate $\sqrt{\tau}$ is indeed maximal.

For $\beta=\tau ,$ the two successors of a point in the overlapping set $D$ will either both return after 3 generations, as for $x=\frac12 ,$ or will return after 2 and $\ge 4$ generations, as in Figure \ref{gm}, or will return later. The case with 2 and 4 gives the highest growth, and Figure \ref{gm} shows that this case can be iterated with itself forever, so this yields the maximum growth rate.  Many other finite orbits combine the 3-3-return with the 2-4-return so that the growth rate is between $\sqrt[3]{2}$ and $\sqrt{\tau}.$  \hfill $\Box $\medskip

{\it Remarks.} Note that the $F_k$ increase exponentially.  Already for $p=15,$ we get $F_9=34$ cycles while the five-fold concatenation of two intersecting 3-cycles would provide only $2^5=32$ cycles. For $p=29$ we obtain 987 cycles on a basic set of 57 points.  

We extended by $11$ if there were more words ending with zero, and by 00 if there were more words ending with 1, in order to obtain the maximum number of cycles. We can also make a random choice of 00 or 11 in each induction step. This will generate other finite orbits with the same period and same number of points, but with a smaller number of cycles, for instance $m=4$ for $p=7,$ and thus with a smaller growth rate. 

For even periods,  we get cycle mixtures without root, as for $p=4$ in Figure \ref{gm}. They are not related to exact overlap of intervals. \medskip

On the right  of Figure \ref{gm} we show the preperiodic orbit of $x=t_2^2+t_2^4.$ Both image points of $x$ eventually land in the fixed points 0 or 1. For each point in this orbit, the number of successors in generation $q$ grows linearly with $q,$ so the growth rate assumes its minimal value 1. The local dimension becomes $\frac{\log 2}{\log \tau}\approx 1.4404.$ Such examples are discussed in Sections 7-9. The minimum 0.94 and maximum 1.44 of local dimension for the Fibonacci measure are well-known, cf. \cite{F4,Hu}. Let us note how they determine the appearance of $\nu$ in our figures.  \medskip

{\it Remark on interpretation of local dimension.}\\ 
By definition \eqref{dimdef}, a measure $\mu$ has local dimension $d$ at a point $x_0$ if the intervals $[x-s,x+s]$ have measure $\approx C\cdot s^d$ for small $s.$ Here $C$ is a positive constant, or a slowly varying function of $s$ like $-\log s.$ A numerical approximation looks like $C\cdot (x-x_0)^{d-1}$ near $x_0.$  For $d>1,$ this function is zero at $x_0.$ For $d=1.5$ the function is like $C\sqrt{x}$ at $x_0=0.$ In Figures \ref{B5} and \ref{B6} such points $x_0$ are seen as dark blue spots. 

For $d<1,$ we have a pole of the function, a red spot in our figures. For $d=0.94,$ however, the singularity looks quite harmless, like $C/\sqrt[17]{x}$ at zero. We can assume $C\approx 1.$ Values of this function are larger than 10 only if $|x|<10^{-17}$ which is below usual numerical thresholds. This explains why our histograms like Figure \ref{drei} do assume rather small values even for Pisot parameters. Singular Bernoulli convolutions are not characterized by large values! It is rather the large variation of the approximating functions, numerically expressed as $\sum_i |\phi (x_{i+1})-\phi(x_i)|$, which indicates the singularity.

A similar remark holds for zeros. If $d$ is only slightly larger than 1, the zero will be hardly visible. For $d=1.06,$ the maximum local dimension $\frac{\log 2}{\log \beta}$ at $t= 0.52 ,$ the function near $x_0$ will look like $\sqrt[17]{x}$ at zero, where values smaller than 0.1 are possible only for $|x|<10^{-17}.$ For parameters near $t=\frac12 ,$ all Bernoulli convolutions appear as strictly positive densities even though it was proved in \cite{JSS} that there is a Cantor set of zeros, with dimension converging to 1 for $t\to\frac12 .$ See Figures \ref{tribo} and \ref{spitze} below.

\section{Network orbits require weak Perron slope}
In general, network orbits in a branching dynamical system are exceptions. The condition that $x$ will repeat in its orbit is an equation of the form $g_{w_n}\cdot ...\cdot g_{w_1}(x)=x.$  Two equations of this type for the same point $x$ will severely restrict the choice of mappings.  This will be made precise for the case of Bernoulli convolutions. The first part of our statement is related to Sidorov's Theorem 2.1 and Proposition 2.6 in \cite{Si9} and Theorem 3.6 in \cite{Si3}. The
role of Perron numbers was emphasized by Thurston  \cite{Th}. He  proved that weak Perron numbers $\beta$ are in one-to-one correspondence $\beta=\exp (h)$ with topological entropies $h$ of maps from $[0,1]$ to itself with finitely many local maxima and minima which all have finite orbits. The second part of our theorem was indicated by the remarks after Theorem 1.2 in \cite{Th}.

\begin{Theorem}[Condition for existence of network orbits] \label{nonex}\hfill\\
Suppose that in a Bernoulli convolution with slope $\beta$ there is a point $y$ such that the orbit of $y$ contains two different paths from $y$ to cycles or finite invariant sets. Then $\beta$ must be an algebraic integer.  If the orbit of $y$ is finite then $\beta$ is a weak Perron number.
\end{Theorem}

{\it Proof that $\beta$ is algebraic. }    
Without loss of generality we assume that the two paths differ at the first step. Thus $y$ belongs to the overlap region, 
and both $g_0(y)$ and $g_1(y)$  lead to a cycle or periodic component.
If $y$ itself is not periodic, there are 0-1-words $v,w,v',w'$ with  $v_1\not=v'_1$ and
\[ g_w(x)=x\quad\mbox{ for }\quad x=g_v(y)\qquad\mbox{ and }\qquad g_{w'}(x')=x'\quad\mbox{ for }\quad  x'=g_{v'}(y)\ .\]
Without loss of generality we assume that the length $m$ of the word $v$ is not smaller than the length $m'$ of $v'.$
When $y$ itself is periodic we have the special case that either $v'$ is the empty word, $y=x'$ and $w'_1\not= v_1,$ or both $v$ and $v'$ are empty, $y=x=x'$ and $w'_1\not= w_1.$ We shall call this the single period and double period case, respectively.

For $w=w_1...w_n\in\{ 0,1\}^n$ we have
$g_w(x)=g_{w_n}\cdot ...\cdot g_{w_1}(x)=\beta^nx +(1-\beta) \sum_{i=1}^{n}\beta^{n-i}w_i \ .$
We solve $g_w(x)=x$  for $x$ and $g_{w'}(x')=x'$ with $w'=w'_1...w'_{n'}$ for $x'.$
\[ x=\frac{\beta -1}{\beta^n-1} \sum_{i=1}^{n}\beta^{n-i}w_i\quad\mbox{ and }\quad x'=\frac{\beta -1}{\beta^{n'}-1} \sum_{i=1}^{n'}\beta^{n'-i}w'_i\ .\]
The equation $x=g_v(y)=\beta^my-(\beta -1)\sum_{j=1}^m \beta^{m-j} v_j$ with $v=v_1...v_m$ and $x'=g_{v'}(y)$ with $v'=v'_1...v'_{m'}$ are both solved for $y$:
\begin{equation} \frac{1}{\beta^m}(x + (\beta -1)\sum_{j=1}^m \beta^{m-j} v_j)\quad =\quad y\quad =\quad  \frac{1}{\beta^{m'}}(x' + (\beta -1)\sum_{j=1}^{m'} \beta^{m'-j} v'_j)\ .\label{m0}\end{equation} 
We multiply both sides with $\beta^m,$ substitute $x$ and $x'$ and cancel $\beta -1$ :
\begin{equation*}
\frac{1}{\beta^n-1}\cdot \sum_{i=1}^{n}\beta^{n-i}w_i +\sum_{j=1}^m \beta^{m-j} v_j = 
\beta^{m-m'}\left( \frac{1}{\beta^{n'}-1}\cdot \sum_{i=1}^{n'}\beta^{n'-i}w'_i  +\sum_{j=1}^{m'} \beta^{m-j} v'_j  
\right) \ .
\end{equation*} 
Multiplying with $\beta^n-1$ and $\beta^{n'}-1$ we obtain polynomials with integer coefficients.
\begin{eqnarray}
(\beta^{n'}-1)\left( \sum_{i=1}^{n}\beta^{n-i}w_i + (\beta^n-1)\sum_{j=1}^m \beta^{m-j} v_j \right)=\notag \\
(\beta^n-1)\beta^{m-m'}\left(  \sum_{i=1}^{n'}\beta^{n'-i}w'_i  + (\beta^{n'}-1)\sum_{j=1}^{m'} \beta^{m-j} v'_j\right)  \ . 
\label{schnitt}\end{eqnarray}
With $N=n+n'-1,$ the leading terms are 
\[ \beta^Nw_1-\beta^{N+m}v_1\  \mbox{ on the left side and }\quad \beta^{N+m-m'}w'_1-\beta^{N+m}v'_1 \ \mbox{ on the right side. }\]
If $m'\ge 1, $ then $m\ge 1$ because of $m\ge m' ,$ and the leading exponent $N+m$ is greater $N$ and $N+m-m'.$ The coefficients $v_1$ and $v'_1$ of $\beta^{N+m}$ on both sides differ by 1 according to the assumption of the theorem, so $\beta$ is the root of a polynomial with integer coefficients and leading coefficient 1.\vspace{1ex}

If $m'=0,$ we have $y=x'$ and the sum with $v'_j$ in \eqref{m0} disappears. For  $m>0,$ the leading term $\beta^{N+m}$ has coefficients $v_1$ and $w'_1$ which now must be different, zero and one. For  $m=0$ the $v_i$ also disappear and we are left with coefficients $w_1$ and $w'_1$ for the leading term $\beta^N,$ which in this case differ since they correspond to application of $g_0$ and $g_1$ to $y$ in the first step. In all cases the leading coefficient is $\pm 1,$ so $\beta$ is an algebraic integer.   \hfill $\Box $\vspace{2ex}

{\it Remark. }Not all choices of words $v,w,v',w'$ are possible in Bernoulli convolutions since equation \eqref{schnitt} will not always have a solution in $(1,2).$ This will be clarified below.  Moreover, without further assumptions on $y$ we cannot expect that $\beta$ will become a weak Perron number. Already the double period case with $n=n'$ and arbitrary $w,w'$ leads to all polynomials with coefficients in $\{ -1,0,1\} ,$ and $p(x)=x^7-x^5+x^3-x^2-1$ has a real root at 1.166 and complex roots with modulus 1.24. \vspace{2ex}

{\it Proof that $\beta$ is weak Perron. } Here we need not require that the whole orbit of $y$ is finite. Let $h$ be the function on $[0,1]$ with $h(x)=g_0(x)$ for $x<y$ and $h(x)=g_1(x)$ for $x>y.$ Only $h(y)=\{ g_0(y),g_1(y)\} $ is two-valued. \emph{We assume that the orbit of $y$ under $h$ is finite.}\\
Since $g_0(z)>z$ and $g_1(z)<z$ for all $z,$ the orbit can be written in ordered form as
$g_1(y)=x_0<x_1<...<x_n=g_0(y),$ and $y=x_k$ with $0<k<n.$ We let $J_i=[x_{i-1},x_i]$ for $i=1,...,n$ and define a matrix \[ M=(m_{ij})\qquad\mbox{  with }\quad m_{ij}=1\quad\mbox{ if } J_j\subseteq h(J_i)\qquad\mbox{ and }\quad m_{ij}=0\quad\mbox{ otherwise. } \] 

This is the adjacency matrix of a graph where an edge leads from $i$ to $j$ whenever $J_j\subseteq h(J_i).$ Since at least one edge starts from each $i,$ there is a vertex set $W\subseteq\{ 1,...,n\}$ such that no edge leads from $W$ to other vertices, and for any two vertices $i,j\in W$ there is a directed path from $i$ to $j.$ (Write $i\preceq j$ if $i=j$ or there is a directed path from $i$ to $j,$ and let $i\sim j$ if $i\preceq j$ and $j\preceq i.$ Take a maximal element $i^*$ with respect to $\preceq$ and let $W$ be the equivalence class of $i^*.$) 

So the adjacency matrix $M_W=(m_{ij})_{i,j\in W}$ of the subgraph induced by $W$ is irreducible. Since the length of $h(J_i)$ is $\beta$ times the length of $J_i,$ the vector of lengths of the $J_i$ with $i\in W$ is an eigenvector of $M_W$ with eigenvalue $\beta .$ By the theorem of Perron and Frobenius this implies that $\beta$ is a weak Perron number. \hfill $\Box $\vspace{2ex}

{\it Remark. }One can show $W=\{ 1,...,n\}$ but this is not needed here. We do not know whether $\beta$ must be a Perron number when $y$ has finite orbit under $\{ g_0,g_1\}.$  Under the weaker assumption that $y$ has finite orbit under $h,$ this will not hold, as shown by the example of $y=\frac12$ for $\beta =\sqrt{2}.$

\section{Points with unique addresses}
Now we shall analyse the parametric family of BCs introduced in Sections 1 and 2. First we consider points $x$ which have a unique address $\pi^{-1}(x)$ in the fractal construction of $\nu_t .$ In other words, for each level $k,$ the point $x$ belongs to a single piece $I_w$ of level $k.$ In terms of our branching dynamical system $G=\{ g_0,g_1\} ,$ there is only one successor of $x$ for each level $k.$ That is, repeated images of $x$ by $g_0$ and $g_1$ will never hit the  overlap interval $D=[t,1-t].$ Typical examples are the points $x=0$ with address $000...=\overline{0}$ and $x=1$ with address $\overline{1},$ for any $t<1.$\smallskip

{\bf The curves for period 2.}
Erd\"os, Jo\'o, and Komornik noted around 1990 that there are no other points with unique address when $t$ is larger or equal the golden mean parameter $t_2.$  The overlap interval {\bf D} and its images ${\bf D}_{0^k}, {\bf D}_{1^k}$ with $k=1,2,...$ will cover the whole interval $(0,1).$ As a consequence, each point $x\in (0,1)$ has a continuum of addresses \cite{EJK, Si3}.  
For $t< t_2 ,$ however, the points
\begin{equation} \textstyle
y=\frac{t}{1+t}<t \mbox{ and } z=1-y=\frac{1}{1+t}>1-t \quad\mbox{ fulfil }\quad g_0(y)=z \mbox{ and } g_1(z)=y\, .
\label{period2}\end{equation}  
So their orbits avoid $D,$ and they have unique addresses $\overline{01}$ and $\overline{10} .$ They form a periodic orbit of length 2 for $t<t_2.$ For $t=t_2 ,$ the points $y$ and $z$ are on the boundary of $D$ and have a countable number of addresses, as is easy to check.

If a point $x$ has a finite or countable number of addresses, the growth rate of its successors is $\rho=0$ and by Proposition 4, the local dimension of $\nu$ at $x$ equals $\frac{\log 2}{\log\beta}>1 .$ A direct consequence is that the density $\lim_{\epsilon\to 0} \frac{\nu [x-\epsilon, x+\epsilon]}{2\epsilon}$ of $\nu$ at $x$ exists and is 0. According to the remark at the end of Section 5, points with unique addresses {\it indicate themselves} by dark blue color in our figures. The curves $y(t)=\frac{t}{1+t}$ and $z(t)=1-y(t)$ for $\frac12\le t\le t_2$ are apparent in Figure \ref{B5}, in particular for $t$ near $t_2$ where they end at the intersection of the dashed line with the black lines. There are several other blue curves which end near the dashed line, for instance $t\cdot y(t).$ \smallskip

{\bf Curves for period 4.} There are many other points with unique address
when $t$ becomes smaller! The overlap region will then be smaller and other periodic points appear outside $D.$ As Figure \ref{stufe4} indicates, the next parameter where a point with unique address separates from {\bf D} is around 0.57. At $t=s_2\approx .5698,$ the inverse of the Pisot number $\phi_2\approx 1.7549$ with minimal polynomial $\beta^3-2\beta^2+\beta-1$ (Figure \ref{bc}). It is easy to check that $\beta^5=\beta^4+1$ and hence $g_{1001}(t)=t.$  For $t<s_2$ we have the following periodic orbit of length 4 outside the overlap region $D.$
\begin{equation} 
\frac{t^2}{1+t^2}\ ,\ \frac{t}{1+t^2}\ ,\ 1-\frac{t}{1+t^2}\ ,\ \frac{1}{1+t^2} \label{period4}
\end{equation}
The unique addresses of these points are $\overline{0011},\,  \overline{0110},\, \overline{1001},\,  
 \overline{1100}\, .$ For $t \in [\frac12 , s_2),$ formula \eqref{period4} defines four curves of points with unique address, density zero and maximal local dimension  $\frac{\log 2}{\log\beta}$  which indicate themselves in our figures by their dark blue color. There are many image curves under contraction maps $f_w$ which in Figure \ref{B5} are blue curves which end at $t=s_2.$\smallskip

{\bf The basic doubling scenario.}
For $s_2<t<t_2,$ however, no further points with unique address exist because for these parameters, the horn ${\bf D}_{(01)^k}$ with $k=1,2,...$ will cover the area between the curve $y(t)$ and the main horn {\bf D} in Figures \ref{stufe4},\ref{B6}. The parameters $t\le t_2$ where new points with unique addresses appear form a decreasing sequence connected with the addresses $\overline{01},\overline{0110}, \overline{01101001}$, $\overline{0110100110010110}.$ The addresses converge to the well-known Thue-Morse sequence, and the $t$-values converge to a parameter $t_{KL}\approx ,5595$ found by Komornik and Loreti \cite{KL} and shown to be transcendental by Allouche and Cosnard \cite{AC}.  

{\bf The sets $A_t.$} For each $t>t_{KL}$ the set $A_t$ of points with unique addresses is countable, for $t\le t_{KL}$ it is uncountable. Starting with basic contributions of Daroczy and Katai \cite{DaK}, Glendinning and Sidorov \cite{GS} and Kall\'os \cite{Kal} many authors have studied the topological structure of $A_t$ and the Hausdorff dimension of $A_t.$ The work of de Vries and Komornik \cite{dVK} collects results up to 2009. More recent papers include Sidorov \cite{Si9}, Jordan, Shmerkin and Solomyak \cite{JSS}, Kong and Li \cite{KoLi}, and Baker \cite{Baker}. The study of sets $A_t$ was inspired by work on univoque numbers \cite{DaKa,AF,KL,dVK} which will be discussed below. \medskip

{\bf Connection with one-dimensional dynamics.}
The doubling scenario with the Thue-Morse sequence is known from the dynamics of one-dimensional unimodal maps, where Milnor and Thurston used a bit different type of address  \cite{MT,CE}, and from the binary coding of exterior angles of the Mandelbrot set, cf. \cite{PR}, where  exactly the addressing is the same as here. The words $01, 0110,...$ describe bubbles of the Mandelbrot set which hit the real axis. Here they describe horns which intersect the main horn {\bf D}. A related result of Allouche, Clarke and Sidorov \cite{ACS} says that the Sharkovskii ordering of periodic points in one-dimensional maps agrees with the ordering of periodic points outside the overlap region which arise for decreasing $t.$  Another analogy with continued fractions was recently found by Tiozzo and co-authors, see \cite{Ti}.

It is easy to explain the reason for the connection between Bernoulli convolutions and one-dimensional real and complex maps.
The doubling map $g:[0,1)\to [0,1)$ is defined by
\[ \textstyle g(x)= 2x \mbox{ mod }1\ ,\quad\mbox{ that is, }\  g(x)=2x \mbox{ for }0\le x< \frac12\quad\mbox{ and }\
g(x)=2x-1 \mbox{ for }\frac12\le x< 1 \ .\]
In complex dynamics it is proved that the map $h(z)=z^2+c$ induces the map $g$ on the external rays of the Julia set for each $c,$ and thus on the external rays of the Mandelbrot set. This was a basic fact for the studies of quadratic maps by Douady, Hubbard and Thurston in the 1980s, cf. \cite{PR}. 

In our context, the following basic statement says that $g$ is conjugate to the dynamics of $G=\{ g_0,g_1\}$ for each $t<t_2$ and for all points $x$ which have their image outside $D.$ In other words, $g$ describes the dynamics of $G$ where it is single-valued, for each possible $t.$  \vspace{10ex}

\begin{Proposition}[Conjugacy of Bernoulli and doubling map]\hfill\\ 
Let $t\in (\frac12,1),$ and let $F_t$ denote the cumulative distribution function of $\nu_t,$ which means $F_t(x)=\nu_t[0,x]$ for $0\le x\le 1.$
Then $F_t$ defines a conjugacy between the action of $G$ on $[0,1]\setminus D$ and the doubling map $g$ on a corresponding subset of $[0,1].$ That is, 
\[ F_t\cdot G(x)=g\cdot F_t(x)\qquad\mbox{ for }\quad x\in [0,1]\setminus D\, .\]
\label{conjlemma}\end{Proposition}

{\it Proof. } We use the special case \eqref{bern1} of the definition of $\nu_t.$ For $x<1-t$ we get
\[ F_tg_0(x)=\nu_t[0,g_0(x)]=\nu_t(g_0[0,x])=2\nu_t[0,x]=gF_t(x)\ .\]
For $x>t$ we have $G(x)=g_1(x)$ and $ F_tg_1(x)$ is transformed as follows;
\[ \nu_t[0,g_1(x)]=1-\nu_t[g_1(x),1]=
1-\nu_t(g_1[x,1])=1-2\nu_t[x,1]=2F_t(x)-1=gF_t(x)\ .
\quad \Box \]

In order to get a conjugacy between dynamical systems, we have to restrict ourselves to points $x$ for which the orbit under $G$ does not intersect $D.$ This means that $x$ has unique address and implies $t\le t_2\approx 0.618.$  Let us recall the Milnor-Thurston concepts of itinerary and kneading sequence, for our case of the doubling map.  Following seminal work by Parry 1960, these 01-sequences were widely used in connection with $\beta$-expansions, see for example  \cite{EJK, GS, AF, dVK}. Binary numbers have the advantage that eventually periodic sequences are identified with rational numbers, and the simple arithmetic of the doubling map $g$ can be used, as it is done for the external angles of the Mandelbrot set \cite{PR}. $g$ acts as shift map on binary representations of numbers in $[0,1].$

\begin{Definition}[Binary itineraries, kneading sequences, address curves]\hfill\\ 
A 01-sequence $b_1b_2...$ and the corresponding binary number $b=.b_1b_2...=\sum_{k=1}^\infty b_k2^{-k}$ are called
itinerary if the orbit closure of $b$ under the doubling map does not contain the point $\frac12 .$  If $b$ itself realizes the minimal distance, that is, no number  $g^{(k)}(b)=.b_kb_{k+1}...$ with $k=1,2,...$ is nearer  to $\frac12$ than $b,$ we call $b$ a kneading sequence. For each 01-sequence $b,$ the function 
\begin{equation}
x_b(t)=\frac{1-t}{t}\cdot \sum_{k=1}^\infty b_kt^k \label{kneaf} \end{equation}
is called the address curve corresponding to $b$ (compare \eqref{fwgw},\eqref{pipi} in Section 2).
\end{Definition}

Binary
itineraries are exactly those 01-sequences which do have neither $\overline{0}$ nor $\overline{1}$ in their orbit closure under the shift (recall that $\frac12$ has the binary adresses $1\overline{0}$ and $0\overline{1}$). In other words, these sequences do not contain $n$ consecutive equal symbols 0 or 1, for some integer $n.$ The kneading sequence corresponding to a given itinerary $b$ can be obtained by determining the orbit closure of $.b_1b_2...$ under the doubling map, and taking the point (or one of the two points) nearest to $\frac12 .$ If $b$ itself is a kneading sequence, then so is $1-b.$ So it suffices to study the case $b<\frac12 ,$ that is, $b_1=0.$ 

Here are a few examples. The sequence $01^101^201^301^4...$ is not an itinerary since the sup of points in the orbit less than $\frac12$ is $\frac12 .$ Itineraries form a dense set, however: for any given 01-word $w,$ the sequence $w010101...$ is an itinerary. In general, the minimal distance of the orbit to $\frac12$ is not realized by a point of the orbit itself, as for $001(01)001(01)^2001(01)^3...$ where the sup of points less than $\frac12$ is $\frac13 .$  Thus we have to take the orbit closure. Kneading sequences, which realize the optimal distance, are rather rare: $\frac13=.\overline{01}$ is the smallest one. \smallskip

{\bf The parameter $t^*.$} For a kneading sequence $b$ with $b_1=0$ we can easily determine the parameter $t^*$ at which $x_b(t)$ enters the overlap region {\bf D}. We solve $x_b(t^*)=1-t^*.$ Thus $t^*=t^*(b)$ is the solution of the equation
\begin{equation}
\sum_{k=1}^\infty b_{k+1}t^k=1\ .
\label{tstar}\end{equation}
There is exactly one solution in $[0,1]$ since the left-hand side is increasing, with value 0 at $t=0$ and value at least one for $t=1.$ When the curve $x_b(t)$ has entered {\bf D} it will remain there: it cannot cross the upper border because the left-hand side does not exceed $\frac{t}{1-t}.$  For a kneading sequence $c$ with $c_1=1$ these calculations would be slightly more complicated. So we use symmetry and say that $t^*(c)=t^*(1-c)$ is the solution of \eqref{tstar} with $b_k=1-c_k.$

Now consider an itinerary $c$ which is not a kneading sequence. Then $t^*(c)$ is defined as the smallest value $t>\frac12$ for which the $G$-orbit of $x_c(t)$ intersects the overlap interval $D.$ The address curves of different points of this orbit cannot intersect for $t<t^*(c)$ since this would contradict Proposition \ref{conjlemma}.  Thus, since all address curves are continuous, even smooth, the address curve of the kneading sequence $b$ of $c$ will be the first member of the orbit which enters {\bf D}, and we calculate $t^*(c)=t^*(b)$ from the kneading sequence $b.$ 

For any itinerary, $t^*(c)$ marks the right endpoint of the dark curve $x_c(t)$ in Figures \ref{B5} and \ref{B6}. At the same parameter where the curve $x_b(t)$ of the kneading sequence $b=.b_1b_2...$ enters {\bf D}, the address curves of a corresponding itinerary $c=.w_1...w_nb_1b_2...$ will enter the horn ${\bf D}_{w_1..w_n}.$ At this point the properties of all these curves will change: as soon as the orbit under $G$ intersects {\bf D}, we have multivalued dynamics.\smallskip

Binary notation assigns rational numbers to eventually periodic 01-sequences, a method used to denote external angles of the Mandelbrot set \cite{PR}.  The curves $y(t)$ and $z(t)$ in \eqref{period2}, for instance, coincide with $x_{1/3}(t)$ and  $x_{2/3}(t)$ defined above. Other rational functions of the form $x_b(t)$ are given in \eqref{period4}, for $b=\frac15 , \frac25,  \frac35, \frac45$ as is easy to check.

In the context of real unimodal maps, like the quadratic family $q_r(x)=rx(1-x),$ itineraries are the addresses of points, and kneading sequences are those addresses which for some parameter $r$ belong to the critical point. A 01-sequence  addresses  different points in different maps, but in a well-behaved parametric family it can address the critical point only once, and then it disappears \cite{MT,CE}. We adapted the Milnor-Thurston notation to Bernoulli convolutions since the situation is similar. Itineraries are unique addresses of certain points, describing the $b$-quantiles as explained below. At the point $t^*$ they become critical which means that the kneading sequence addresses a boundary point of {\bf D}, and the address curves of other itineraries enter other horns. At $t^*$ the points cease to have a unique address and maximal local dimension. As we shall see, the curves $x_b(t)$ remain important beyond $t^*$ although in Figures \ref{B5} and \ref{B6}, $x_b(t^*)$ seems to be their endpoint.

Milnor and Thurston \cite{MT} also introduced kneading functions
similar to our address curves in order to determine the topological entropy of unimodal maps. The standardizing factor $\frac{1-t}{t}$ in \eqref{kneaf} comes from our choice of mappings $g_0,g_1$ which define all measures on $[0,1].$  The following theorem says that the address curves of itineraries are just the parallel dark curves in Figures \ref{B5} and  \ref{B6}. \vspace{10ex}

\begin{Theorem}(Points with unique address, kneading functions and quantile curves)\label{quant}\hfill 
\begin{enumerate}
\item[ (i)] For each $t\in (\frac12, t_2),$ the set $A_t$ of points with unique address agrees with all values $x_c(t)$ outside {\bf D}. More exactly,
\begin{eqnarray*}  A_t&=&\{ x_c(t)\, |\, c \mbox{ itinerary and } t<t^*(c)\, ,\\ 
&& \mbox{ or } t=t^*(c)\, ,\mbox{ and the $g$-orbit of $c$ does not contain the kneading sequence} \}\ .\end{eqnarray*}
\item[(ii)]  For each itinerary $b=.b_1b_2...,$ the kneading function $x_b(t)$ represents the $b$-quantile of all Bernoulli measures $\nu_t$ with $t\le t^*(b)$ :
\[\textstyle F_t(x_b(t))= \nu_t[0, x_b(t)] = b \quad\mbox{ for }\ {\textstyle\frac12}\le t\le t^*(b)\ .\] 
\end{enumerate}
\end{Theorem}

{\it Proof of (i). } We fix some $t<t_2$ and consider the branching dynamical system $G$ for this $t.$ Take any itinerary $c$ with $t<t^*(c).$ By the definition of $t^*$ the $G$-orbit of $x_c(t)$ does not hit the overlap interval D. This also holds in the other case $t=t^*(c)$ by Proposition \ref{conjlemma} since then the kneading sequences correspond to the endpoints of $D.$ Thus $x_c(t)$ has a unique address and belongs to $A_t.$ \\ 
To show the reverse inclusion, take a point $x\in A_t,$ and let $c_1c_2...$ denote its unique address, explained in Section 6. Thus there is only one $G$-successor of $x$ in each generation: $x_1=g_{c_1}(x), x_2=g_{c_2}(x_1), ...$ In other words, the $G$-orbit of $x$ does not hit $D=[1-t,t].$ By Proposition \ref{conjlemma}, the orbit of $c=.c_1c_2...$ under the doubling map does not hit the interval $[s,1-s]$ where $s=\nu_t[0,1-t]<\frac12 .$ Thus $c$ is an itinerary and $x=x_c(t).$ Moreover, $t\le t^*(c)$ since the $G$-orbit of $x$ does not hit $D.$\\ 
If the $g$-orbit of the binary number $c$ contains a corresponding kneading sequence $b$ then by Proposition  \ref{conjlemma} and the definition of $t^*$ we have $t< t^*(c).$ This holds in particular for all eventually periodic sequences $c,$ the case which is our main interest here. However, if the $g$-orbit of $c$ does not contain a corresponding kneading sequence, $t=t^*(c)$ is possible. \hfill $\Box$\smallskip

{\it Example.}  Let $c=(01)(0110)(01101001)...$ be the concatenation of approximating words of the Thue-Morse sequence. Then $t^*=t_{KL}$ is the Komornik-Loreti point, and $b$ is the Thue-Morse sequence which for this parameter is an address of the endpoint  $1-t^*$ of $D$ which has $1\overline{0}$ as second address. However $x_c(t^*)$ has only the address $c.$ Taking $c=(01)^{k_1} (0110)^{k_2} (01101001)^{k_3}...$  with nonnegative integers $k_i$ we obtain uncountably many sequences with the same properties. \smallskip

 Our address curves are analytic functions for $0<t<1.$ The assertion that points with unique address come on smooth curves seems to be new. We consider address curves as structural elements of the Bernoulli scenario which provide a smooth pattern at least outside {\bf D}. 
A consequence of (i) is that  
\[ A_t\subseteq \overline{A_t}\subseteq \{ x_c(t)\, |\, c \mbox{ itinerary and } t\le t^*(c)\} =:\tilde{A_t}\ ,\]
and all three sets coincide if the $G$-orbits of $t$ and $1-t$ hit the interior of $D,$ as for instance for $t\in (s_2,t_2).$ If the orbit of $1-t$ avoids the interior of $D$ (and so has either two or countably many addresses), the sets differ by a countable set of preimages $f_w(t), f_w(1-t)$ outside $D.$  In this case $1-t$ can be an accumulation point of $A_t,$ as in the above example. Then $\overline{A_t}=\tilde{A_t}$ and the difference occurs between the closure $\overline{A_t}$ and $A_t.$ When $1-t$ is not an accumulation point, as for $t=t_2,$ then $A_t$ is closed and $\tilde{A_t}\setminus A_t$ is countable.
The differences often do not matter, for instance in calculating Hausdorff dimension or measures.

\section{Geometry of itineraries and kneading sequences}\label{geom}
The topological structure of $A_t$ is the subject of a comprehensive paper by de Vries and Komornik \cite{dVK} based on properties of greedy and quasi-greedy $\beta$-expansions. It collects earlier results of several authors and considers the general case $\beta >1.$ In \cite{dVK,GS,Si9}, ${\cal U}_q, {\cal V}_q$ denote the non-standardized version of $A_t, \tilde{A_t},$ and $q$ means $\beta .$

The goal of this section is to give a direct self-contained, geometric and dynamic view of the $A_t,$ for our case $1<\beta\le 2,$ to prove (ii) of Theorem \ref{quant} and indicate connections with one-dimensional dynamics. Roughly speaking, all combinatorial and topological questions boil down to the study of the angle-doubling map. For $m-1<\beta\le m$ with $m>2$ we would have to take the map $mx$ mod 1.

To each set $A_t$ we can assign the set $I_t$ of initial points $x_c(\frac12 )$ of the curves $x_c$ with $x_c(t)\in A_t.$ Assertion (ii) says that $F_t$ maps $A_t$ to $I_t.$ Note that $F_t$ is a homeomorphism on $[0,1]$ since the measure $\nu_t$ has no point masses, and each interval has positive measure. This suggests that the study of the \emph{sets $ \tilde{A_t}$ can be replaced by the study of sets $S_b$ of binary itineraries} subordinated to a given kneading sequence $b.$ The $S_b$ have a rather simple structure. In the following statement $f_0(x)=x/2, f_1(x)=(x+1)/2,$ and $f_w$ for a 01-word $w$ is the corresponding composition of mappings, as explained in Section 6.

\begin{Proposition}[Structure of itineraries subordinated to a kneading sequence]\label{St}\hfill\\ 
Let $b<\frac12$ be a kneading sequence and $J=(b, 1-b).$ Let $S_b$ be the set of all itineraries for which the orbit under the doubling map does not meet $J.$ Then there is a set $W$ of 01-words such that the intervals $J_w=f_w(J)$ with $w\in W$ are pairwise disjoint, and
\[  [0,1]\setminus S_b =\bigcup_{w\in W} J_w\ . \]
Actually, $W$ consists of all words $w=w_1...w_m$ such that for $k=1,...,m-1$ the binary number $.w_k...w_m$ is not in the interval $[b,1-b).$ Moreover, $b$ is either an isolated point of $S_b,$  or an accumulation point of the points $f_w(b)$ with $w\in W.$ The first case happens if and only if $b$ is a periodic binary number of the form $.\overline{v_1...v_n(1-v_1)...(1-v_n)}.$
\end{Proposition}

{\it Proof. }  
First we show that the intervals $J_v=f_v(J)$ for all 01-words form a nested sequence: any two of them are either disjoint, or one is contained in the other. Namely, if $J_v$ would contain an endpoint of $J,$ for some word $v$ of length $n>0,$ then $g^n(b)$ or $g^n(1-b)$ would be in the interior of $J$ which cannot happen since $b$ and  $1-b$ are kneading sequences. 
Next, if $J_u$ contains an endpoint of $J_v,$ but not the whole set $J_v$ for some words $u=u_1...u_m$ and $v=v_1...v_n$ then we must have $u_1=v_1$ since $f_0(0,1)\cap f_1(0,1)=\emptyset ,$ and $J_{u_2...u_m}$ and $J_{v_2...v_n}$ satisfy the same overlap relation as $J_u$ and $J_v.$ Repeating this argument we come to the case where one of the intervals is $J,$ which was shown to be impossible. 

Thus the $J_v$ are nested, and we let $W$ be the set of words $w$ for which $J_w$ is not contained in any other $J_v.$ Note that an itinerary $c\in S_b$ cannot be in any $J_v,$ and any point outside the union of all $J_v$ must be an itinerary. Thus $\bigcup_{w\in W} J_w$ is the complement of $S_b.$ The inductive argument above says that $J_u\subset J_v$ for $u=u_1...u_m$ and $v=v_1...v_n$ if and only if $m>n$ and $u_i=v_i$ for $i=1,...,n$ and $J_{u_{n+1}...u_m}\subset J.$ The characterization of $w\in W$ says that this cannot happen. 

$b$ is isolated if it is the right endpoint of an interval $J_w$ for some nonempty word $w\in W.$ Thus $b=f_w(1-b),$ and by symmetry $1-b=f_{1-w}(b),$ consequently $b=f_wf_{1-w}(b).$ Otherwise, every interval $(b-\epsilon,b)$ with $\epsilon >0$ contains an interval $J_w$ with $w\in W,$ and $f_w(b)$ is the left endpoint of that interval. \hfill $\Box$\bigskip

{\it Remarks.} (Structure and Hausdorff dimension of $A_t$)
\begin{enumerate}
\item The condition for $W$ allows us to determine the number $a_m$ of holes $J_w$ in the Cantor set $S_b$ on every level $m,$ and the growth rate $\rho=\lim_{m\to\infty} \frac{\log a_m}{m}.$ The Hausdorff dimension of $S_b$ then is $\log\rho/\log 2,$ cf. Kong and Li \cite{KoLi}.  
\item The above proof works for the sets $A_{t^*}$ with $t^*=t^*(b)$ when we replace $J$ by $(1-t^*,t^*).$ It is also possible to relate $S_b$ to $A_{t^*}$ by assigning to each starting point $(\frac12 , c)$ of a blue curve with $c\in S_b$  the terminal point $(t^*,x_c(t^*))=(t^*,1-t^*)$ with $t^*=t^*(b).$ Thus $A_{t^*}$ has the same similarity structure  with mappings of contraction factor $t$ as $S_b$ with contraction factor $\frac12 .$ As a consequence \cite[Theorem 2.6]{KoLi}, for $t=t^*$
\begin{equation}\label{dim}
\textstyle {\rm dim}\, A_t=\frac{\log\rho}{-\log t}=\frac{\log\rho}{\log \beta} \end{equation}
\item In the case where $b$ is isolated, $\rho$ can be determined explicitly, similar to Proposition \ref{grow}, cf. \cite{DaK,Kal}. In particular, when $W$ consists of all words which do not have $k$ consecutive equal symbols, for some $k\ge 3,$ then $\rho=t_{k-1}.$ This holds for $b=(2^{k-1}-1)/(2^k-1)$ and $t^*(b)=t_k,$ and was used in \cite{JSS,GS} to show ${\rm dim}\, A_t\to 1$ for $t\to\frac12 .$
\item The $S_b$ form an increasing sequence, so their dimension increases with $b.$ However, Kong and Li \cite{KoLi} observed that ${\rm dim}\, A_t$ \emph{sometimes decreases} with decreasing $t.$ When $b$ is periodic, there is an interval $(b,b')$ which does not contain kneading sequences, see Proposition \ref{dub} below. For $t\in (t^*(b'),t^*(b))$ the above proof, with $(1-t,t)$ instead of $J,$ still yields equation \eqref{dim}, and $A_t$ has the same similarity structure while $\log\beta$ increases with decreasing $t.$ 
On such intervals, $ {\rm dim}\, A_t$ smoothly decreases with decreasing $t.$ 
\end{enumerate}

Let us now study
the structure of the set $K$ of all kneading sequences $b<\frac12$. In Figure \ref{B6},  kneading sequences are represented as endpoints of blue curves on the line $y=1-t.$ To any starting point $(\frac12 , b)$ for a kneading sequence $b<\frac12$ we assign the `terminal' point $(t^*,x_b(t^*))=(t^*,1-t^*)$ with $t^*=t^*(b).$ However, many blue curves correspond to itineraries which are not kneading sequences. The set $A_{t^*}$ is the intersection of the vertical line $t=t^*$ with all blue curves.  Going back on the curves, we get the corresponding set $S_b$ on the vertical line $t=\frac12 .$ 

\begin{figure}[h]
\begin{center}
\includegraphics[width=0.8\textwidth]{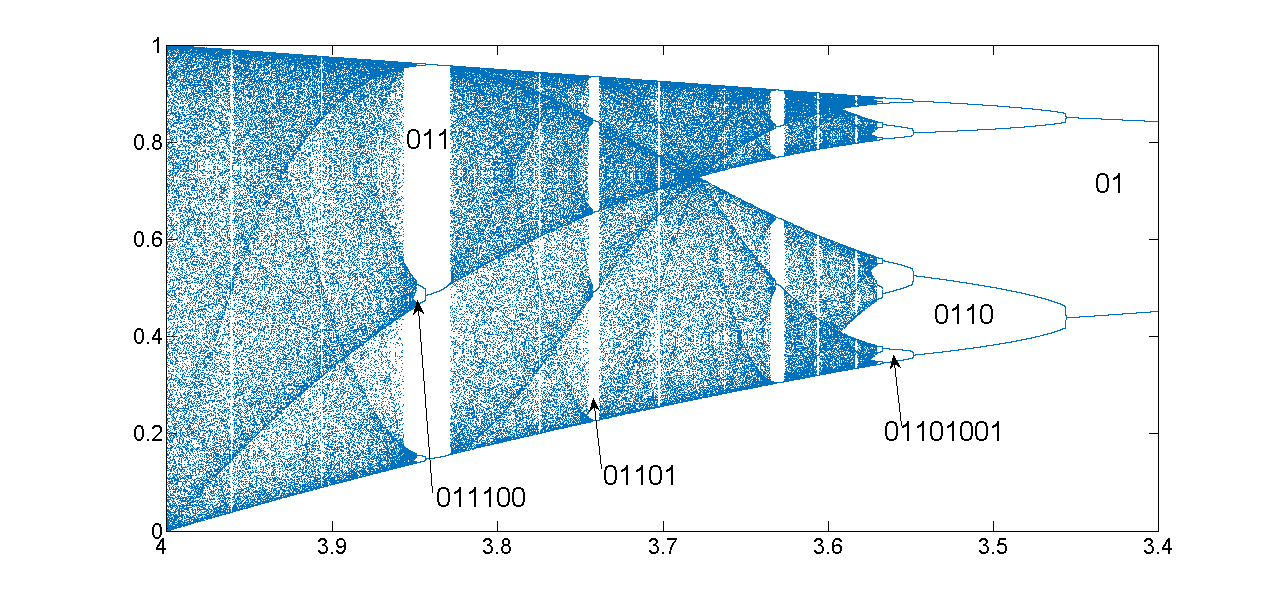}\\
\includegraphics[width=0.8\textwidth]{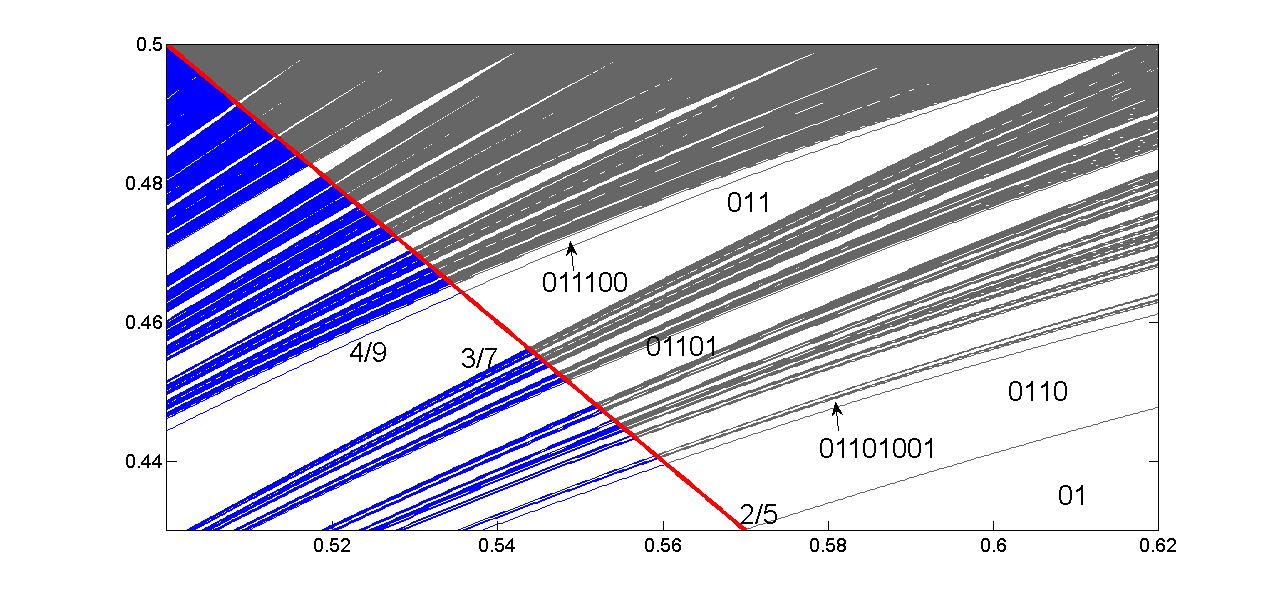}
\end{center}
\caption{Bifurcation diagram of the quadratic family $q(x)=rx(1-x), r\in [3.4, 4]$ with binary itineraries for 'periodic windows' \cite{CE,MT,PR}, and address curves of kneading sequences in our scenario. Combinatorial structures coincide while metric structures differ a lot. In the lower panel, intersection with the line $y=1-t$ represents kneading sequences, and univoque numbers as a subset.}
\label{knea}
\end{figure}

\begin{figure}[h]
\includegraphics[width=0.999\textwidth]{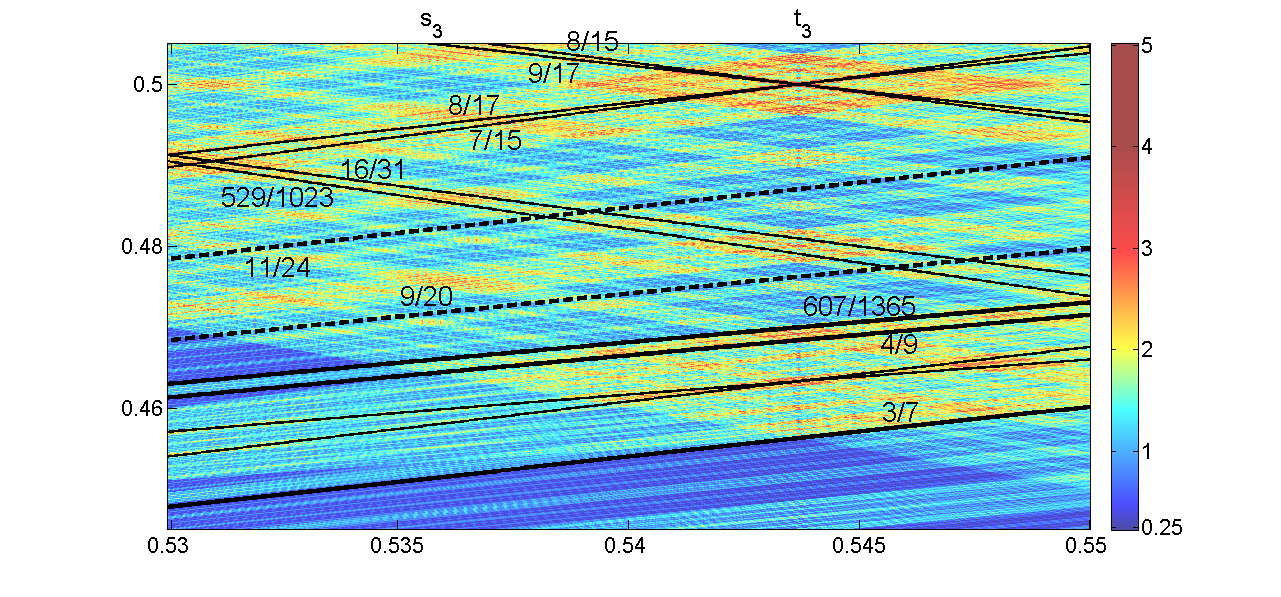}
\caption{A magnification of Figure \ref{B6} near the tribonacci parameter $t_3.$ Address curves of some eventually periodic itineraries are shown in black and extended into the overlap region.}
\label{tribo}
\end{figure}

Figure \ref{knea} shows only the address curves of $b\in K$, except for $b=\frac13$ which is outside range. Each kneading sequence is represented by the intersection of its address curve with $y=1-t,$ and curves are extended up to the line $y=\frac12$ where some of them meet. 
The points $(t^*,1-t^*),$ as well as the kneading sequences $b$ on the left axes, form a Cantor set plus a sequence of isolated points in each cutout interval which always converges to the left endpoint of the interval, which is a Komornik-Loreti point. This phenomenon is well-known from one-dimensional dynamics, where the `bifurcation points' form a geometric sequence studied by Feigenbaum and others, and so are better visible due to the slowdown of the dynamics when the parameter reaches a Feigenbaum point. In Figure \ref{knea} we use binary words to denote the periodic windows of the Feigenbaum diagram, which corresponds to the notation of the Mandelbrot set \cite{PR}, not to the original Milnor-Thurston addressing. With this notation, the correspondence between parameter windows with stable periodic orbit, and gaps in the Cantor set of kneading sequences is one-to-one. In the lower picture, however, the distance of isolated points to the left endpoint of the gap interval decreases with order $t^{2^n},$ that is, with double exponential speed, so that only the very first isolated kneading sequences can be distinguished from the corresponding Komornik-Loreti limit. Compare Figure \ref{B6} for the basic case $.\overline{01},\ .\overline{0110}, ...$ with $t^*=t_2,s_2,...\to t_{KL},$ and the thick black lines in Figure \ref{tribo} for  $.\overline{011}=3/7$ with $t^*=t_3,$ \ $.\overline{011\,100}=4/9$ with $t^*=s_3,$ and $.\overline{011\,100\ 100\,011}=607/1365.$ \medskip

{\it Remark.} (Potential of correspondence between quadratic maps and $\beta$-expansions)\hfill\\ 
There is a vast literature on the dynamics of real unimodal maps. Some results carry over to $\beta$-expansions and Bernoulli convolutions. One instance is Sharkovskii ordering of orbit lengths of real continuous maps which was proved to hold for the dynamics of $G$ on the $A_t$ by Allouche, Clarke, and Sidorov \cite{ACS}. 
The existence of absolutely continuous invariant measures is a central problem in one-dimensional dynamics, and Misiurewicz parameters play a similar part as Garsia numbers here. Sometimes the setting of linear expansive maps is simpler than quadratic maps: a point of period $n$ is a root of a polynomial of degree $n$ while for quadratic maps the degree would be $2^n.$ \medskip

The following well-known fact explains the period-doubling sequences of isolated kneading addresses. It says that every periodic sequence $b\in K$ is isolated from the left in $K,$ and if it has the form $b=.\overline{w(1-w)}$ for some word $w,$ it is also isolated from the right. See \cite[Section 3]{AF} for related facts references.

\begin{Proposition}[Period-doubling lemma]\label{dub}\hfill\\
Let $b=.\overline{0b_2...b_n}$ be a periodic kneading sequence. There is no kneading sequence between $b$ and $b'=.\overline{0b_2...b_n1(1-b_2)...(1-b_n)}.$ \end{Proposition}

{\it Proof. } Let $c>b$ be a kneading sequence which we try to make as small as possible. Let $k$ be the first index for which $c_k>b_k.$ For $k>n+1$ the sequence $c$ would not be kneading, so we take $k=n+1, c_{n+1}=1.$  Since $.c_{n+1}c_{n+2}...\ge 1-c=.1(1-b_2)...(1-b_n)0...$ we inductively conclude $c_{n+j}\ge 1-c_j$ for $j=2,3...$ and take $c_{n+j}=1-c_j$ as optimal choice. This yields $c=b'.$ Since $b$ is a kneading sequence, $b'$ also fulfils this condition.  
\hfill $\Box$\medskip

A number $\beta\in (1,2)$ is called \emph{univoque} if the equation
\begin{equation}\label{univo}
\sum_{k=1}^\infty c_k\beta^{-k} =1 \quad\mbox{ with coefficients }\quad c_k\in\{0,1\}
\end{equation}
has only one solution. This concept was thoroughly studied in connection with $\beta$-expansions \cite{EJK,DaKa,DaK,AF,dVK}. A basic fact is that \emph{the set of univoque representations agrees with the set of non-periodic kneading sequences} \cite[Theorem 1]{EJK}, \cite[Theorem 2.3]{AF}. Thus univoque parameters $t=1/\beta$ are obtained by removing all isolated points and right endpoints of gaps from $K,$ considered as set of $t^*(b).$ This implies the topological statements in \cite{dVK}, where the sets are denoted  ${\cal U}, {\cal V}.$ In particular, the set of univoque numbers is not compact. 

\begin{Proposition}[Univoque numbers and kneading sequences \cite{EJK}, Theorem 1]\label{univoque}\hfill\\
Let $\beta\in (1,2)$ and $c_1c_2...$ a 01-sequence which fulfil  \eqref{univo}.  Then $\beta$ is univoque if and only if
$b=.0c_1c_2...$  is a non-periodic kneading sequence. If this is true, $\beta =1/t^*(b),$ that is, $\beta$ is determined by \eqref{univo}.   \end{Proposition}

{\it Proof. } If $\beta$ is univoque with coefficients $c_k,$ and $b=.0c_1c_2...$ then $x_b(t)=1-t$ by definition.
$b$ is a kneading sequence and $t=t^*(b)$ since the $G$-orbit of $1-t$ does not meet $D,$ by the univoque condition. The assumption of periodic $b=.\overline{0c_1...c_k}$ implies 
\[ 1-t=x_b(t)=\frac{1-t}{1-t^{k+1}} \sum_{j=1}^k c_jt^j \quad \mbox{ and }\quad  t^{k+1}+\sum_{j=1}^k c_jt^j=1\  ,\]
which contradicts the univoque assumption. So $b$ is non-periodic.\\
Conversely, for a non-periodic kneading sequence $b<\frac12$ and $t=t^*(b)$ the $G$-orbit of $1-t$ will never return to $D.$ So $1-t$ has only two addresses $b=.0b_2b_3...$ and $.1\overline{0}.$ This implies that $b_2,b_3,...$ is the only possible coefficient sequence for the representation of 1, and $\beta$ is univoque.\hfill $\Box$\medskip

Univoque Pisot numbers were studied by Allouche, Frougny and Hare \cite{AF} for which Theorem \ref{finiteorbits} says that $b$ must be eventually periodic. Thus a Pisot number is univoque if and only if $b$ is preperiodic.  An extensive, partly computer-assisted search in \cite{AF} proved that the preperiodic case is quite rare for Pisot numbers:

\begin{Theorem}[Smallest univoque Pisot numbers \cite{AF}]\label{univoquepisot}\hfill\\
The smallest accumulation point of univoque Pisot numbers is the root $\chi\approx 1.9052, t\approx .5256$ of $x^4-x^3-2x^2+1.$ This is the only accumulation point of Pisot numbers with preperiodic kneading sequence. Only two univoque Pisot numbers are below $\chi ,$ at $1.8800$ and $1.8868\, .$ 
\hfill $\Box$ \end{Theorem}
 
The address of $\chi$ is $.011\overline{10}=\frac{11}{24}.$ It will show up at Theorem \ref{central} below.
The following approximation of kneading sequences by periodic ones will be needed in the proof below.

\begin{Proposition}[Approximation of kneading sequences by periodic ones]\label{doub}\hfill\\ 
Let $b=.b_1b_2...<\frac12$ be a nonperiodic kneading sequence and $J=(b, 1-b).$ The orbit of $s_n=.\overline{b_1...b_n}$
under the doubling map does not meet $J$ if $b_{n+1}=1$ and there is no $k<n$ such that $(1-b_k)...(1-b_{n+1})=b_1...b_{n-k+2}.$ There are infinitely many $n$ with this property.
\end{Proposition}

{\it Proof. } The assumption $b_{n+1}=1$ implies $s_n<b\, ,$ and $g^m(s_n)<g^m(b)\le b$ whenever $b_{m+1}=0$ and $1\le m\le n.$ Thus $g^m(s_n)$ cannot belong to $(b,\frac12 ]$ for any $m\ge 1.$  Now suppose that $g^m(s_n)\in (\frac12 , 1-b)$ for some $m<n.$ Then $b_{m+1}=1$ and $.b_{m+1}...b_n0 < 1-b$ while  $.b_{m+1}...b_nb_{n+1}...\ge 1-b$ holds because $b$ is a kneading sequence. Hence $.(1-b_{m+1})...(1-b_n)1>b$ and $.(1-b_{m+1})...(1-b_n)0\le b .$  So this last word must be a prefix of $b$ which contradicts the last assumption. The first statement is proved.  

Induction will show that the assumption is fulfilled for infinitely many $n.$ Suppose that at some stage, we have problems to find the next $n$ since there is a word $b_kb_{k+1}...$ which is the complement of a prefix of $b.$ The word cannot extend to infinity because then $b$ would be periodic, with period $2(k-1).$ Take the smallest $m$ such that 
$(1-b_k)...(1-b_m)$ is a prefix of $b$ but $1-b_{m+1}\not= b_{m+2-k}.$ By the kneading property, $b_{m+1}=1=b_{m+2-k}.$ There can be no $k'>k$ such that $(1-b_{k'})...(1-b_m)0$ is a prefix of $b$ since $(1-b_{k'})...(1-b_m)1$ is a subword of $b$ whenever $k\le k'<m.$ We found the next $n=m$ and can continue.
\hfill $\Box$\medskip

{\it Proof of Theorem \ref{quant}, (ii). } Consider first an eventually periodic 01-sequence $b_1b_2...$ It is an itinerary if it does not end with $\overline{0}$ or $\overline{1}.$ The corresponding binary number $b=.b_1b_2...$ is rational but not of the form $k/2^n.$ So the orbit under $g$ is finite and Proposition \ref{mark} applies. We obtain all values $F_{\frac12}(g^k(b))$ from the equations of the Markov partition. Of course we know that $F_{\frac12}(x)=x$ since $\nu_\frac12$ is the uniform distribution. However, Proposition \ref{conjlemma} says that for $t<t^*(b)$ the orbit of $x_b(t)$ is also finite, and gives rise to exactly the same equations for the Markov partition. In particular $F_t(x_b(t))=F_\frac12 (b)=b .$ It is easy to see that this is also true for $t=t^*(b)$ since Bernoulli convolutions have no point masses.

Next, we prove (ii) for nonperiodic kneading sequences $b.$ Proposition \ref{doub} says that the binary number $b=.b_1b_2...$ can be approximated by periodic sequences $s_n$ with $t^*(s_n)\ge t^*(b).$ So the address curves $x_{s_n}(t)$ converge to $x_b(t)$ for $\frac12\le t\le t^*(b).$ Since (ii) holds for the $s_n$ and the $\nu_t$ are nonatomic measures, (ii) must hold for $b.$ 

Now we show (ii) for itineraries. Suppose $b$ is an itinerary with $F_t(x_b(t))=b$ and $t\le t^*(b),$ and  $c=b/2$ is also an itinerary with $t^*(c)=t^*(b).$ Then $x_c(t)=t\cdot x_b(t)$ by definition,  and Proposition \ref{conjlemma} implies $F_t(x_c(t))=b/2=c$ since
\[  2 F_t(x_c(t))=F_t\cdot g_0(x_c(t))= F_t(x_b(t))=b\ .\]
A similar equation holds when $c=(b+1)/2$ is an itinerary. By recursion, (ii) is carried over from a kneading sequence $b$ to all itineraries $d$ with $t^*(d)=t^*(b)$ which have $b$ in their $G$-orbit.  Proposition \ref{St} says that such points $d$
approximate also the other itineraries which have $b$ only in the closure of their $G$-orbit. So (ii) holds for them since $F_t$ is a homeomorphism.  \hfill $\Box$\medskip

The quantile property (ii) directly implies that address curves cannot intersect at $t<t^*,$ that is, outside the horns ${\bf D}_w .$ It gives some information on $\nu_t$ for all $t\le t_2,$ no matter whether singular or absolutely continuous. Moreover, for periodic or preperiodic binary numbers $b$ the address curve is a rational function of $t.$

This also holds for  $b=.b_1b_2...b_m1\overline{0}$ which is not an itinerary. In this case, $x_b(t)$ is the middle curve of the horn ${\bf D}_w$ with $w=b_1...b_m$ which is mapped by $g_w$ onto the middle curve $y=\frac12$ of ${\bf D}.$ As long as these curves enter no horns of lower level, they also represent the $b$-quantile of all $\nu_t.$ For example, the curve $x(t)=\frac{t}{2}$ in the middle of ${\bf D}_0$ represents the $\frac{1}{4}$-quantile for all $t\le\frac{2}{3}$ where $x(t)$ hits the lower border $1-t$ of ${\bf D}.$ For ${\bf D}_{01}$ the equation $\beta^2x+1-\beta =\frac12$ leads to the middle curve $x(t)=t-\frac12 t^2$ which defines the $\frac{3}{8}$-quantile for $t\le 2-\sqrt{2}.$  The general case is proved by induction and \eqref{bern}. In contrast to other quantile curves, the density of $\nu$ on these curves is not zero.

\section{Density zero inside the overlap region} 
We saw that the single-valued dynamics of $G$ outside $D$ boils down to the study of the doubling map.
The multi-valued dynamics inside $D$ is much less understood.  Our view is that the address curves remain structural elements within {\bf D}, and their intersections determine zeros and poles of the two-dimensional density $\Phi .$

We start with a look on Figure \ref{tribo}. Some address curves have been extended beyond the parameter $t^*$ where they enter {\bf D}. As mentioned above, the three thick curves in the lower part correspond to the period 3 doubling scenario starting with $.\overline{011}=3/7\approx .429.$  The period 4 doubling scenario is represented by the addresses
$.\overline{0111}=7/15\approx .467$ and $.\overline{0111\, 1000}=8/17\approx .471$ and by their counterparts
$.\overline{1000}=8/15$ and $.\overline{1000\, 0111}=9/17$ which come from above. Period 5 doubling is represented only by two curves with $b>\frac12 ,$ namely  $.\overline{10000}=16/31$ and  $.\overline{10000\, 01111}=529/1023.$
Each number denotes the starting point of the respective curve at $t=\frac12$ which is not visible in the figure.

The two thin lines without numbers arise as images under $f_{01}$ of the $7/15$ and $8/15$ curve. Their adresses are preperiodic, $.01\overline{0111}=22/60\approx .367$ and $.01\overline{1000}=23/60\approx .383.$ 
Finally, there are two dashed lines with preperiodic addresses $.01\overline{1100}=9/20=.45$
and $.011\overline{10}=11/24\approx .458,$ and these curves pass through regions with small $\Phi$-values, which will be proved rigorously. The periodic address curves, however, pass through high $\Phi$-values, at least at those places were such curves intersect. 

Since the addresses are eventually periodic, the curves are rational functions of $t.$ 
So the $t$-value of an intersection point of two such curves is a root of a polynomial, and $\beta=1/t$ is a Perron number according to Theorem \ref{nonex}. All curves of the period 4 doubling scenario intersect at $(t_3,\frac12 ),$ and the period 3 curves meet at $(t_2,\frac12 ),$ as can be seen in Figure \ref{knea}. A similar fact can be checked for period $n$ and the multinacci parameter $t_n.$ 
Actually, all other periodic kneading sequences, like $.\overline{01101}$ give rise to a doubling scenario with infinitely many address curves which will all meet at a point $(t,\frac12 )$ where $\beta=1/t$ is a Perron number. 

Intersections of periodic address curves $x_b(t), x_c(t)$ will be treated in the next section. Here we consider those intersections which give small $\Phi$ values. When $b$ and $c$ have a common prefix, say $b=ab'$ and $c=ac'$ with $a=a_1...a_n,$ then $x_{b'}(t)$ and $x_{c'}(t)$ will intersect at the same parameter(s) $t$ as $x_b(t)$ and $x_c(t),$ and the $x$-value of the latter intersection point is in the $G$-orbit of the former.
For this reason we shall assume $b_1\not= c_1.$

\begin{Proposition}[Intersection points with two addresses]\label{intersec0}\hfill\\ 
Let $b,c$ be nonperiodic kneading sequences with $b_1=0$ and $c_1=1,$ let $t\in (\frac12, t_2]$ be the smallest parameter for which $x_b(t)=x_c(t)=y,$ and suppose the $G$-orbit of $y$ does not return to $D.$ Then $y$ has exactly two addresses, the local dimension of $\nu_t$ at $y$ assumes its maximal value $\log 2/\log\beta ,$ and the density of $\nu_t$ at $y$ is zero. All points in $D$ with two addresses are obtained in this way.
\end{Proposition}

{\it Proof. } The assumption implies that $y$ is in $D,$ and both $g_0(y)$ and $g_1(y)$ have unique addresses. The value of the dimension comes from Theorem \ref{dime}.
 \hfill $\Box$\medskip

\begin{figure}[h]
\includegraphics[width=0.999\textwidth]{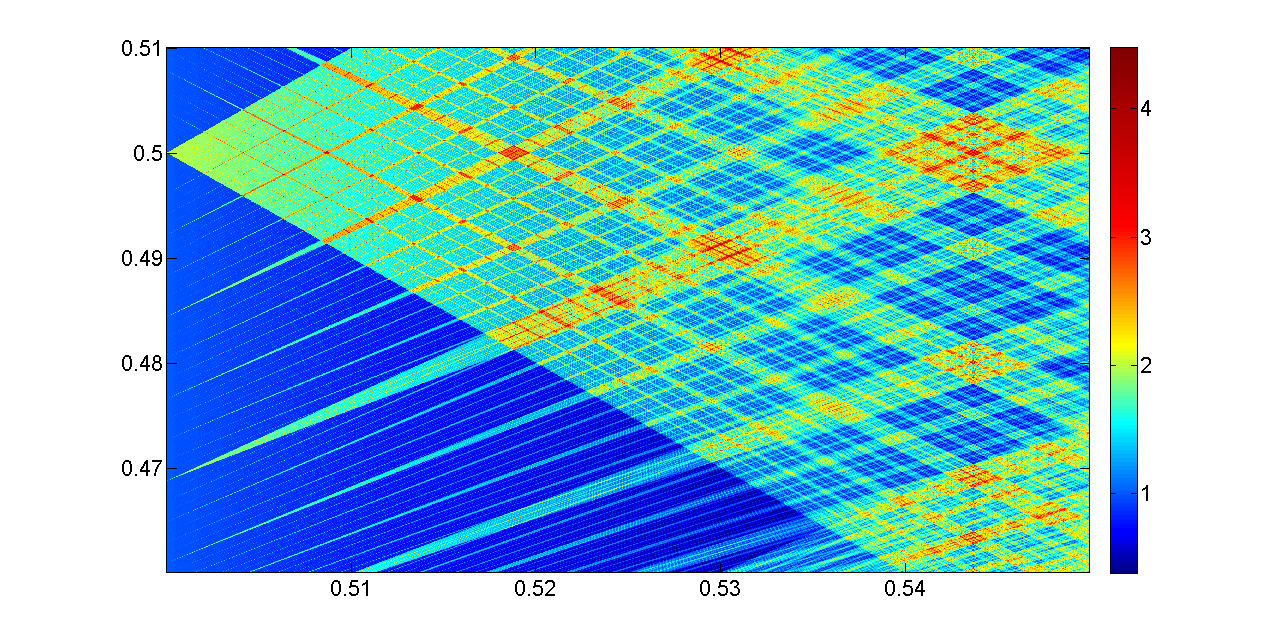}
\caption{Near $t=\frac12$ the structure of $\Phi$ looks self-similar, with lots of zeros within {\bf D}. Landmark parameters are the multinacci numbers $t_n,$ the doubling numbers $s_n,$ and other Perron parameters.}
\label{spitze}
\end{figure}

Figure \ref{spitze} indicates that this case happens very often when $t<t_3.$ Points with two addresses seem to form Cantor carpets within {\bf D}. For $b>0.469,$ the binary sequence starts $b=.01111b_6...$ so that $x_b(t)=(1-t)(t+t^2+t^3+t^4+...)=t-t^5+...$  Thus near $t=\frac12 ,$ the curves $x_b(t)$ become parallel to the line $y=t,$ and the curves $x_c(t)$ become parallel to $y=1-t.$  A better explanation of the apparent product structure of the dark carpets within {\bf D} is an observation of Sidorov \cite[Lemma 2.2]{Si9} which we reformulate for our setting. 

\begin{Proposition}[Condition for points with two addresses]\label{intersec1}\hfill\\ 
For fixed $t<t_2,$ there is a one-to-one correspondence between points $y\in D$ with exactly two addresses and pairs of points $x,x'\in A_t$ with $x-x'=\beta -1 .$ In particular, points with two addresses exist if and only if $\beta -1$ belongs to the difference set $A_t-A_t.$
\end{Proposition}

{\it Proof. } Given $y,$ we let $x=g_0(y)$ and $x'=g_1(y).$ Given $x,x'$ let $y=f_0(x)=f_1(x').$
 \hfill $\Box$\medskip

Since $\beta -1>1-t$ for $0<t<1,$ we must always have $x>\frac12$ and $x'<\frac12 .$ Thus we can construct the $g_1$-image of the points with two addresses as
\begin{equation}\label{intersec}\textstyle 
\mbox{ intersection of }\quad \{ x'\in A_t\, |\, x'<\frac12\}\quad\mbox{ with }\quad  
\{ x+1-\beta \, |\,  x\in A_t,\, x>\frac12\} \ . 
\end{equation}
Using address curves, this can be done for all $t<t_2$ together. This explains the Cantor product structure, but it does not clarify how large the set of points with two addresses is for fixed $t.$ The best answer to this question for $t$ near $\frac12$ was found by Sidorov.

\begin{Theorem}[Points with two addresses for all $t<t_3$ \cite{Si9}, Theorem 4.2]\label{sidorov2}\hfill\\ 
For each $t<t_3,$ there is at least one point with exactly two addresses.
\end{Theorem}

\begin{figure}[h]
\begin{center}
\includegraphics[width=0.7\textwidth]{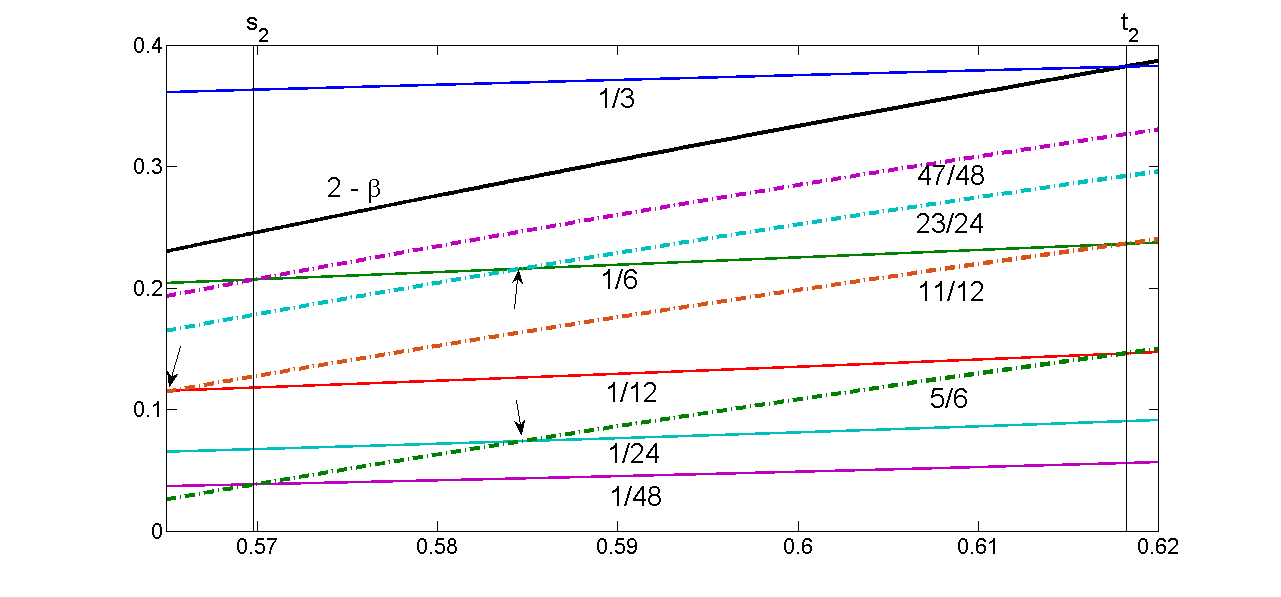}
\caption{Points with two addresses for $0.565\le t\le 0.62$ obtained by intersections of curves $x'(t), x(t)$ given in \eqref{intersec}. The first set is represented by solid lines, the second by dotted ones.}
\label{intersec2low}
\end{center}\end{figure}

{\it Example. }
As Figure \ref{intersec2low} shows, \eqref{intersec} can be easily applied to $s_2<t<t_2$ where $A_t$ is given only by the itineraries  $2^{-k}\frac13=.0^k\overline{01}$ and $1-2^{-k}\frac13=.1^k\overline{10}$ for $k=0,1,2,...$ The shifted curve of $2/3$ is out of range, and only $k\le 4$ is considered. The curves for $k>4$ are between $x_{1/48}$ and 0, and between $x_{47/48}$ and $2-\beta$ so they provide no other intersection points. There are intersection points at $t_2$ where we have no points with unique address, and at the doubling number $s_2.$ The largest parameter $t\approx .5846$ admitting points with two addresses was found by Sidorov \cite[Proposition 2.4]{Si9}.  Figure \ref{intersec2low} demonstrates that this is the only parameter $t>s_2$ with this property \cite[Theorem 2.8]{Si9}. The intersection points are marked by two arrows. 

Figure \ref{intersec2low} shows another intersection parameter $t<s_2$ at the left border, marked by a single arrow. In this case, $x=1-x',$ so for symmetry reasons the point with two addresses must be $y(t)=\frac12 .$ We have  $x'(t)=t^3/(1+t)$ with address $.00\overline{01}$ and $y=f_1(x')$ with address  $b=.100\overline{01}\, .$ Curiously, $1/t^*(b)$ is the Pisot number $\chi$ of Theorem \ref{univoquepisot}. The equation $y(t)=t^4/(1+t) +1-t=\frac12$ with $t\not=1$ yields $2t^2(t+1)=1$ and $t\approx .5652\, ,$ a Garsia parameter:

\begin{Theorem}[Parameters for which $\frac12$ has only two addresses]\label{central}\hfill\\ 
There are uncountably many parameters $t$ for which the central point $\frac12$ has exactly two addresses. The parameter $t=1/\beta \approx .5652 ,$ where $\beta$ denotes the real root  of $x^3-2x -2,$ is the largest accumulation point of this set. The largest isolated parameter is $t\approx .5674 .$
\end{Theorem}

{\it Proof. } Due to symmetry, Proposition \ref{intersec1} simplifies for the central point: $y=\frac12$ admits exactly two addresses if and only if $y+\frac{1-\beta}{2}=1-\frac{1}{2t}=:\varphi (t)$ belongs to $A_t.$ Let $\psi (t)=t^3(1-t)$ denote the lower bound of the horn ${\bf D}_{000}=f_{000}({\bf D})$ with tip at $t=\frac12, y=\frac{1}{16}$ (cf. Figure \ref{stufe4}). We have $\varphi (0.5)=0<\psi (0.5)=1/16$ and $\varphi (0.54)\approx .0741>\psi (0.54)\approx .0724$ so all address curves of itineraries $c<1/16$ with $t^*(c)\ge 0.54$ will intersect the curve $\varphi (t)$ on the interval $[0.5, 0.54].$ There are uncountably many kneading sequences $b<\frac12$ with $t^*(b)\ge 0.54,$ and for each of them $c=000b$ is an appropriate itinerary.

As shown above, $t\approx .5652$ marks the intersection point of $\varphi (t)$ with $x_{00\overline{01}}(t).$
Now $\varphi (t)$ is strictly increasing, and  for $t<s_2\approx .5698$ the address curves of $00(01)^n\overline{0110}$ approximate $x_{00\overline{01}}(t)$ from above, as noted after \eqref{period4}. So $\varphi (t)$ must intersect almost all of these curves. Calculation shows that the largest intersection point below $s_2$ appears for $n=3$ at $t\approx .5674 .$ Other address curves need not be considered since $t^*(.\overline{0110\, 1001})\approx .5604.$
 \hfill $\Box$\medskip

Theorem \ref{central} is somewhat surprising since at first glance one would expect that Bernoulli convolutions cannot have density zero at the center of the interval.  It also answers a question in \cite[Section 5]{Si9}: Let ${\cal B}_m$ denote the set of parameters $t$ which admit points with exactly $m$ addresses. 
Consider the set of ${\cal B}_2\cap [t_{KL}, t_2].$  Is this set discrete? No, it has accumulation points. An argument similar to the proof above shows that $s_2$ is an accumulation point of the set from the left, using address curves of $0^3(01)^n\overline{0011}$ approaching $x_{1/48}=x_{.0^3\overline{01}}$ in the lower left corner of Figure \ref{intersec2low}. Hence $s_2$ is the largest accumulation point of this set.  \medskip

The sets ${\cal B}_m$ for $m>2$ are obviously subsets of ${\cal B}_2.$ They were studied in \cite{Si9,BS} and will not be treated here.  The set  ${\cal B}_{\aleph_0}$ of parameters admitting points with infinite countable number of addresses, however, is not a subset of  ${\cal B}_2.$ Such parameters can be generated similar to Proposition \ref{intersec1}.

\begin{Proposition}[Intersection points with countable number of addresses]\label{intersecaleph}\hfill\\ 
Let $b$ be a nonperiodic and $c$ a periodic kneading sequence with $b_1\not= c_1,$ let $t\in (\frac12, t_2]$ be the smallest parameter for which $x_b(t)=x_c(t)=y,$ and suppose the $G$-orbit of $y$ does not return to $D\setminus\{ y\}.$ Then $y$ has a countably infinite number of addresses, the local dimension of $\nu_t$ at $y$ assumes its maximal value $\log 2/\log\beta ,$ and the density of $\nu_t$ at $y$ is zero. 
\end{Proposition}

{\it Proof. } If $b=b_1b_2...$ and $c=\overline{w}$ with $w=c_1...c_n,$ the addresses of $y$ have the form $w^kb$ with $k=0,1,2...$
 \hfill $\Box$\medskip

We briefly study the interval $t_{KL}<t\le t_2$ where all points of $A_t$ are eventually periodic. For such parameters $t,$ points with two addresses arise iff two preperiodic address curves meet. Points with countably many addresses arise when a periodic kneading sequence meets a preperiodic one, or, perhaps, infinitely many preperiodic curves meet in one point. Theorem \ref{nonex} applies to all these cases, which improves \cite[Theorem 2.1]{Si9}:

\begin{Proposition}[Number of addresses within first doubling scenario]\label{oneorinfinite}\hfill\\ 
If $t_{KL}<t\le t_2$ and $\beta=1/t$ is not a weak Perron number, then each point $y$ has either one address or uncountably many addresses.
\end{Proposition}

In contrast to ${\cal B}_2,$ the set ${\cal B}_{\aleph_0}$ includes $t_2,$ and every other $t$ for which $1-t$ has a periodic address $b=.\overline{0b_2...b_n}$ which is a kneading sequence. These were the non-univoque numbers excluded in Proposition \ref{univoque}.  

We conclude the section with two examples near $t_2.$ Proposition \ref{intersecaleph} applies to the kneading sequences $b=\frac{5}{12}=.01\overline{10}$ and $c=8/15=.\overline{1000}$ , with intersection parameter $t\approx .5951$ and $y\approx .463\, .$ For $c'=135/255=.\overline{1000\, 0111}$ we also get an intersection point with the address curve $x_b(t)$ at $(t',y')\approx (.6045, .4668).$ The $G$-orbit of $y'$ returns to $1-y'$ and is countable although Proposition \ref{intersecaleph} does not apply.

\section{Singularities inside the overlap region} \label{super}
So far we studied points $y$ which have very few addresses - a finite or countable number, while almost all points in $[0,1]$ have an uncountable number of addresses \cite[Theorem 3.6]{Si3}.  These points have an interesting structure but they do not cause $\nu_t$ to become singular.  Now we consider points $y$ which have network-like orbits under $G,$ where at least two cycles in the network have a common point. 
The growth rate of such points is always positive, and they can have an extraordinary number of addresses:   

\begin{Theorem}[Intersections of periodic address curves]\label{intersecper}\hfill\\ 
Let $b=.\overline{0b_2b_3...b_m}$ and $c=.\overline{1c_2c_3...c_n}$ be periodic itineraries, and let the two curves $x_b(t), x_c(t)$ intersect in the point $(s,z)$ inside the overlap region $D.$ Then
\begin{enumerate}
\item[(i)] 
Infinitely many periodic address curves meet in $(s,z).$
\item[(ii)]
The growth rate of the orbit of $z$ is at least as large as the positive solution of 
\begin{equation}\label{mn}x^{-m}+x^{-n} =1\ .
\end{equation}
\item[(iii)]
If the growth rate of the orbit of $z$ exceeds $2s$ then $\nu_s$ does not have a bounded density. A possible density function of $\nu_s$ must be unbounded on each interval, and discontinuous at every point of $[0,1].$
\end{enumerate}
\end{Theorem}

{\it Proof.  (i) } If $z$ is the fixed point of $f_v$ and $f_w,$ with $v=0b_2b_3...b_m$ and   $w=1c_2c_3...c_n,$ then $z$ is also the fixed point of $f_{vw}$ and $f_{wv},$ of $f_{vwwv}$ and $f_{wvvw}$, of $f_{vwwvwvvw}$ and so on. Thus $z$ belongs to the address curves of $.\overline{vw}, .\overline{vwwv}$ etc.  \smallskip

{\it (ii) } We show that the growth rate of a graph containing two directed cycles of length $m$ and $n$ which meet in a vertex $x$ is at least as large as the positive solution of \eqref{mn}. When $z$ is the only intersection point of the cycles and there are no further edges, the number $\gamma^q=\gamma^q(z)$ of successors of $z$ in generation $q$ is given by the equation $\gamma^q=\gamma^{q-m}+\gamma^{q-n}.$ The standard ansatz $\gamma^q=x^q$ leads to \eqref{mn}, and the solution is unique since we look for the spectral radius of a graph and the Perron eigenvalue is unique. When the cycles meet in other vertices and/or there are more edges in the graph, the number of successors of $z$ can only become larger. \smallskip

{\it (iii) } If $\rho$ denotes the growth rate of $z,$ Theorem \ref{dime} says that the local dimension 
is $d=d_z(\nu)=\log (2/\rho) /\log \beta .$ Thus $d<1$ if and only if $\rho >2/\beta =2t.$ In case $d<1$ a possible density function $\phi$ of $\nu$ cannot be bounded in a neighborhood of $z.$ Since $d_{f_w(z)}(\nu)\le d_z(\nu )$ for any 01-word $w,$ and the points $f_w(z)$ are dense in $[0,1],$ the function $\phi$ must be unbounded on each interval. This implies that $\phi$ is discontinuous at each $x\in [0,1].$    \hfill $\Box$\vspace{10ex}

\begin{figure}[h]
\centerline{\includegraphics[width=0.5\textwidth]{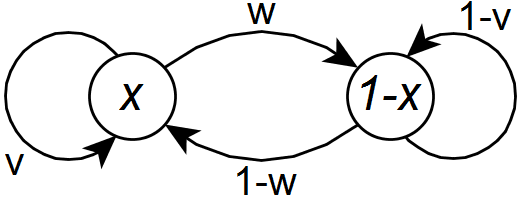}}
\caption{The two-cycle graph for a combined orbit of $x$ and $1-x.$ Edges here stand for a whole path in Section 4, and the two paths to $x$ may share a terminal part, as in Figure \ref{gm}.}
\label{last1}
\end{figure}

\begin{figure}[h]
\centerline{\includegraphics[width=0.999\textwidth]{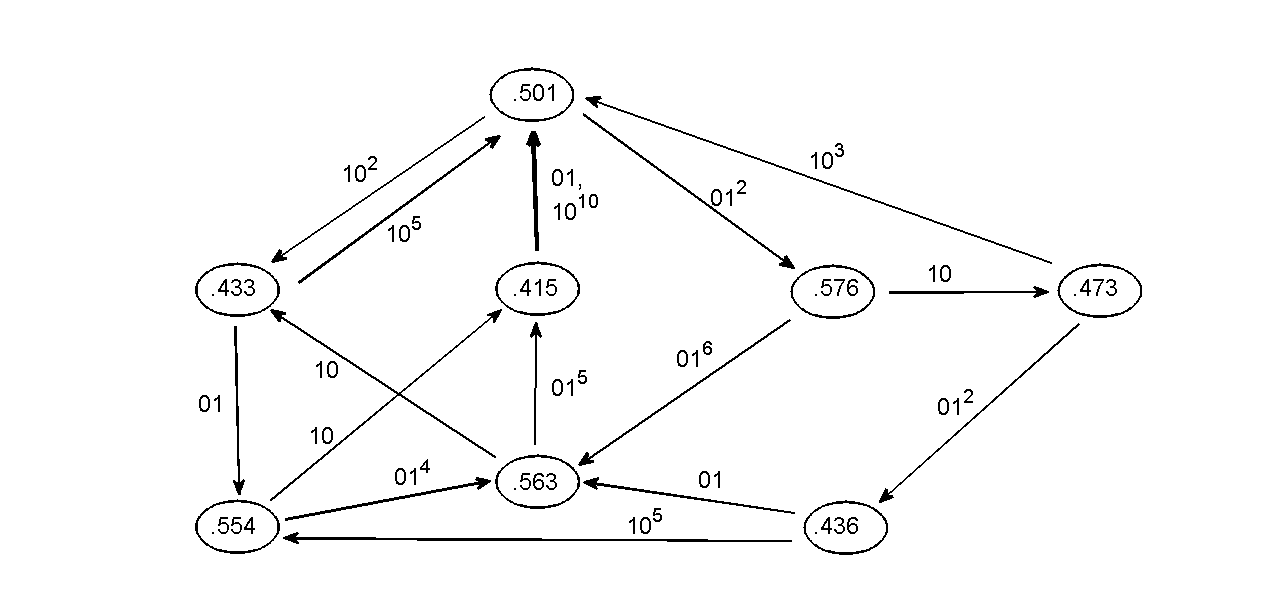}}
\caption{Network orbit with eight branch points for the Pisot parameter $\alpha\approx 1.7049.$ Arrows indicate chains of maps $g_0,g_1.$ Two incoming chains can have a terminal part in common.}
\label{last2}
\end{figure}

{\it Remarks } (Corollaries and improvements of Theorem \ref{intersecper}) \hfill
\begin{enumerate}
\item  For kneading sequences with $c=1-b$ we get in (i) the doubling scenario of Section \ref{geom}. In this case $z=\frac12 ,$ and Figure \ref{knea} shows how the bundle of address curves, representing a gap of the set of kneading sequences, converges to a point $(s,\frac12 )$ where the gap closes. Such intersection points are dense on the line $y=\frac12 .$ In the theory of Milnor and Thurston \cite{MT}, all unimodal maps with addresses from a doubling scenario have the same topological entropy, namely $-\log s.$ This is also the case here for certain $\beta$-expansions. The formulation with words $v,w$ is an immediate generalization, and can be related to topological entropy of Lorenz maps, which will not be done here. In the general case, the addresses need not be kneading sequences. 

\item Here is a little corollary without formal proof. Imagine that $1/s$ is a Pisot number, so that intersection points $z$ of periodic addresses are dense in $D.$ Then the bundles of sufficiently close intersection points $z_1,z_2$ must intersect each other, either left or right of $s.$ Thus a Pisot parameter is always an accumulation point of other parameters where periodic addresses meet, both from the right and from the left. Since there are only few accumulation points of Pisot numbers, most of the approximating parameters will not be Pisot, only weak Perron by Theorem \ref{nonex}.  

\item Combining (ii) with (iii), we can say that there is no bounded density if 
\begin{equation}\label{mns}
(2s)^{-m}+(2s)^{-n}>1 \ . \end{equation}
The case $m=n$ was found by Feng and Wang \cite[Theorem 1.4]{FW} who also gave several conditions under which a density cannot belong to $L^q$ for various $q.$

\item In general, the estimates \eqref{mn} and \eqref{mns} are rather weak. For a given $t_0>\frac12$ there can be only finitely many parameters $s>t_0$ for which these conditions are fulfilled. To see this, assume $m\le n,$ and take $m_0$ so that  $(2t_0)^{-m_0} <\frac12 .$ Then \eqref{mns} can only hold if $m<m_0.$ Moreover, for each given $m<m_0$ there is $n_0$ with  $(2t_0)^{-n_0} <1-(2t_0)^{-m} ,$ which implies that only finitely many $n$ can fulfil \eqref{mns} for $s>t_0$ together with $m.$ Finally, for given $n$ and $m$ there are only finitely many 01-words $v$ and $w$ of length $m$ and $n,$ respectively. For given words $v,w,$ the intersection point is given by a polynomial which has a finite number of solutions.

\item Another type of graph which is more likely to have large growth rate is shown in Figure \ref{last1}.
There is one orbit combining $z$ and $1-z,$ defined by $g_v(z)=z, g_w(z)=1-z,$ due to symmetry. If the lengths of $v,w$ are $m,n,$ respectively, the estimates \eqref{mn} for the growth and  \eqref{mns} for nonexistence of a bounded density remain valid. This is easily shown by replacing $\gamma^q$ with $\gamma^q(z)+\gamma^q(1-z)$ in the proof of (ii). See Example 1 below, and Section 4 where we proved that such an orbit is optimal for the Fibonacci parameter. 
However, the above finiteness argument remains valid also for this double type of graph.

\item For Pisot parameters, there will often be more intricate network-like orbits which lead to growth rates exceeding the estimate \eqref{mn}. See Example 2 below and the network in Figure \ref{last2}. For non-Pisot parameters, we have not found such networks. So it remains possible that beside Pisot parameters, for every $t_0>\frac12$ there are only finitely many numbers $s$ for which there is an intersection point $z$ with growth rate greater $2s.$  
\end{enumerate}

{\it Example 1 } (Perron parameter admitting a very small local dimension)\hfill \\
One way to find examples for Theorem \ref{intersecper} is to prescribe $v,w$ and determine the corresponding $s,z.$ To get large growth rates, we take the double type of graph in Figure \ref{last1} with $v=10000, w=01.$ The equation $g_v(z)=z$ implies $\beta^4(\beta z +1-\beta)=z$ and $z=\frac{\beta^4(\beta-1)}{\beta^5 -1},$ and $g_w(z)=1-z$ yields $z=\frac{\beta}{\beta^2+1}.$
This gives the minimal polynomial $\beta^5-\beta^4-\beta^2 -\beta -1$ with Perron root
$\beta\approx 1.6851$ and two complex pairs of roots of modulus 1.03 and 0.75, respectively. Thus $s=1/\beta\approx .5934.$ The point $z\approx .4389$ has growth rate $\rho \approx 1.2365,$ the real solution of \eqref{mn} for $m=5, n=2.$ By Theorem \ref{dime}, the local dimension of $\nu_s$ at $z$ is $d\approx .9215.$ The minimal local dimension for the Fibonacci parameter is 
$\frac{\log 2}{\log \tau}-\frac12\approx 0.9404,$ cf. \cite{F4,FS} and Section 4, where we have the same graph with $v=1000, w=01.$ The Pisot parameter of the next example gives minimal local dimension 0.9380. Thus our Perron parameter has `larger peaks' than neighboring Pisot parameters! Nevertheless, we do not expect $\nu_s$ to be singular. See Figure \ref{drei}. \smallskip

{\it Example 2 } (Pisot parameter with a typical network orbit)\hfill \\
The Pisot number $\alpha\approx 1.7049$ with minimal polynomial $\alpha^5-\alpha^4-\alpha^3 -1$ was identified by Sidorov as the first parameter for which a point with exactly two addresses exist (Section 9). The number $z\approx .501$ is the fixed point of $g_v$ and $g_w$ with $v=10^210^5$ and $w=01^21010^3.$ Since $m=n=9,$ equation \eqref{mn} gives $\sqrt[9]{2}\approx 1.08$ as lower bound for the growth rate, which is much smaller than $\frac{2}{\alpha}\approx 1.173.$ However, since we have a Pisot number we should have supercritical orbits. It turns out that the orbit of $z$ is a network with eight branching points, and with growth factor $\rho\approx 1.2125.$ This gives the local dimension $0.9380$ which seems minimal for this parameter. The network is shown in Figure \ref{last2}. $\rho$ can be determined from $\det (A-I)=0$ where $A$ is the adjacency matrix of the network, with entries $\rho^{-n}$ instead of 1 for edges $(x,y)$ where $y=g_w(x)$ and $w$ has length $n.$ Networks without short cycles seem to be typical for Pisot numbers with high degree, where conjugates have modulus almost one. For Salem numbers of degree 4 and 6, results of Boyd indicate that the situation could be similar, see Thurston \cite{Th}. Networks can be much larger than Figure \ref{last2}. However, singularities generated by complicated networks seem to have small impact on measures of neighboring parameters! This is a point which requires more study. \smallskip

\begin{figure}[h]
\centerline{\includegraphics[width=0.999\textwidth]{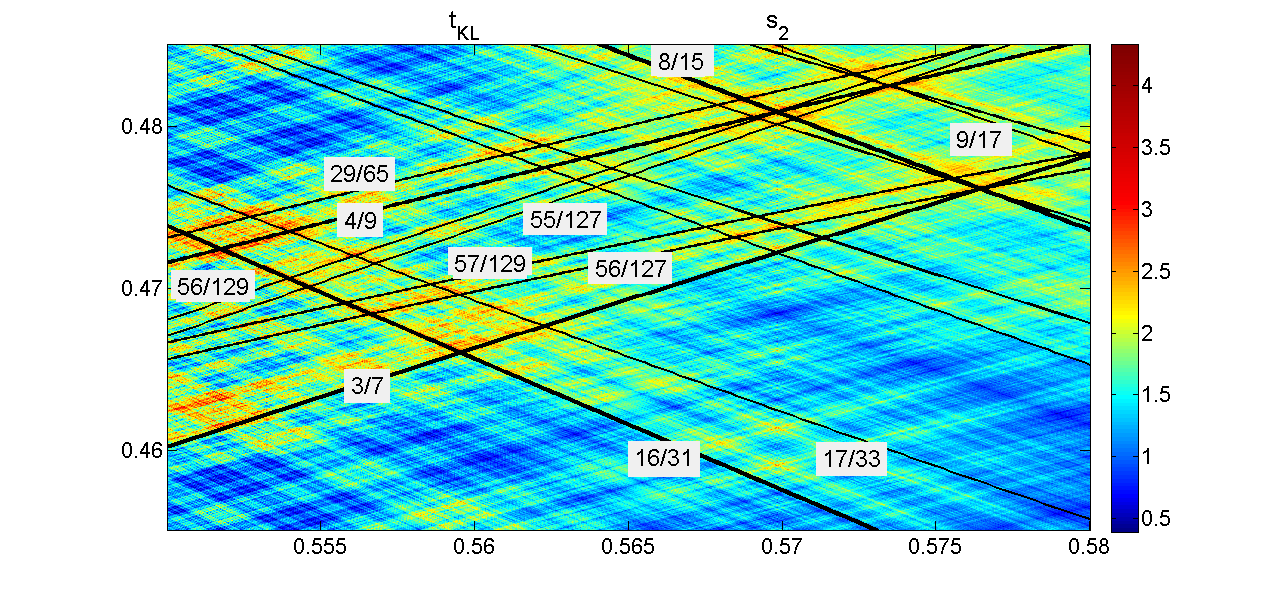}}
\caption{$\Phi(t,x)$ for $t\in [.55,.58],\, x\in [.455,.485]$ with periodic address curves of small period. Apparently all peaks of $\Phi$ are intersections of such curves.}
\label{last3}
\end{figure}

{\it Example 3 } (Generating Perron examples from address curves)\hfill \\
It is possible to determine intersection points by a systematic search of periodic address curves. In Figure \ref{last3} we have chosen a window where the origin of peaks of $\Phi$ is not as obvious as for multinaccci parameters, and have tried to fit the peaks by intersections of periodic address curves. This seems to work! In \cite{ban} there is a picture with less detail.

The window contains the Komornik-Loreti point, and one can distinguish the dark blue sets generated by countably many lines right of $t_{KL}$ from the dark Cantor carpets left of $t_{KL}.$ The Pisot parameter $s_2$ is also in the window, and there are many intersections of address curves at $t=s_2$ as well as near $s_2.$ These are mostly Pisot parameters since $1/s_2$ is a limit point of Pisot numbers (cf. Section 4). The largest peak in the window is at $(s_2, .4809),$ the intersection of curves for $b=4/9=.\overline{011100}$ and $c=8/15=.\overline{1000}.$ This is an orbit of the type shown in Figure \ref{last1}, with $v=1000$ and $w=011.$ The growth rate $\rho\approx 1.221$ is obtained from $\rho^{-3}+\rho^{-4}=1,$ and the local dimension is $d\approx 0.8778,$ much smaller than for the above examples and the Fibonacci number., cf. the appendix in \cite{HHM}.

Another Pisot parameter $s=0.5765$ is determined by the intersection for $b=3/7=.\overline{011}$ and $c=8/15=,\overline{1000}.$ This is just the intersection of a 4-cycle with a 3-cycle, and \eqref{mn} gives the same growth rate as above, which leads to $d\approx 0.8963.$ Figure \ref{last3} illustrates Theorem \ref{intersecper}(i): the address curve of $56/127=.\overline{0111000}$ must pass through the intersection point.

Address curves in Figure \ref{last3} come in bundles, with corresponding binary numbers $3/7, 4/9=28/63, 29/65$ from the left to the upper right, $8/15$ and $9/17$ from top to right, $16/31$ and $17/33$ from upper left to lower right. Rational numbers $k/(2^p-1)$ denote binary representations $.\overline{w_1...w_p}$ of period $p.$ A simple calculation shows that the period-doubling number $.\overline{w_1...w_p(1-w_1)...(1-w_p)}$ is the rational $(k+1)/(2^p+1).$ The next curve of the period-doubling bundle was not drawn since it  almost coincides with the second one.

These are kneading sequences, and their canonical intersection patterns can be seen very clearly in Figures \ref{B6} and \ref{spitze}. Our window is in a parameter region where there are no other kneading sequences of low period. Our hypothesis is that high values of $\Phi$ can be explained by intersection of periodic address curves. To demonstrate this, we take two bundles of itineraries of period $m=7$ which are not kneading: $55/127,\, 56/129$ and $56/127,\, 57/129$ from lower left to upper right. In \cite{ban} we proved that the curves corresponding to $55/127$ and $16/31$ intersect at $s\approx .5546, \, z\approx .4701$ leads to a Perron number of degree 9, and \eqref{mn} with $m=7, n=5$ yields $\rho \ge 1.1237>2t.$ By Theorem \ref{dime}, the local dimension is $d_z(\nu_s)\le .978.$

The curves associated with $56/129$ and $16/31$ intersect at $s'\approx .5540,\, z'\approx .4706.$  The minimal polynomial for $s'$ is $t^8-t^7+t^5+t^4+t^2+t-1.$
The intersection point of any two rational functions is of course the root of a polynomial. Theorem \ref{nonex} says that $\beta =1/s'$ is a weak Perron number. It is not Pisot. We have a graph as in Figure \ref{last1}. \eqref{mn} applies and shows that we have the same lower bound for $\rho$ as for $s.$ Since $s'<s,$ however, we get a slightly better estimate for local dimension: $d_{z'}(\nu_{s'})\le .976.$ Moreover, it turns out that the orbit of $z'$ will visit $D$ once more than expected, and so there will be more successors of $z'$ than estimated by \eqref{mns}. A numerical check showed that the large peaks of $\nu_{s'}$ at the resolution of Figure \ref{drei} can be explained by $z', 1-z'$ and their images under maps $f_v.$

The curves corresponding to $56/127$ and $16/31$ intersect at $s"\approx .5566$ where $1/s"$ is a Perron number of degree 7, and  \eqref{mn} with $m=7, n=5$ again gives the same bound for the growth rate and local dimension $\le .9841$ at $z"\approx .4684.$ Again, the orbit visits $D$ in another point and has further branches including the fixed points of $g_{10^31^20}$ and $g_{10^21010}.$ We have not been able to decide wether this gives rise to a larger growth rate. A numerical approximation shows quite a number of large peaks. Finally, curves of $56/127$ and $17/33$ yield another Perron parameter $s^*\approx .5595$ with a graph like Figure \ref{last1} and again the same growth bound, and still supercritical: $d\le .9928$ at $z^*\approx .4696.$
..
Further Perron parameters are easy to find. Checking for supercritical growth and Pisot/non-Pisot  can be left to computer. Some preliminary work was done, and many non-Pisot parameters found.  Their number rapidly grows when we go with $t$ nearer to $\frac12 .$ We did not find larger network orbits for Perron parameters, however.\medskip

{\it Example 4 } (Transient growth of an orbit)\hfill \\
Consider the Perron number $\beta\approx 1.7924$ with minimal polynomial $x^5-x^4-x^3-x^2+x-1.$ The complex conjugates have modulus around 1.11 and 0.68. No bounded density can exist by Feng and Wang \cite{FW} since $\frac12$ is the fixed point of $g_{100010}.$ Since the orbit of $\frac12$ does not contain other points of $D,$ equation \eqref{mn} gives the growth rate $\rho=\sqrt[6]{2}$ and the local dimension $d\approx .9898$ at $\frac12 .$ However, there are other local dimensions and a non-trivial multifractal spectrum which will not be proved here. At a resolution of 2.5 million bins,  the approximation of $\nu$ has a value 1.7 at $\frac12$ and the maximum value 3.2 at $z\approx .46737.$ It turns out that 
\begin{equation} \label{last}
g_{01110}(z)=z\ \mbox{ and }\ g_{1000}(z)=y\ \mbox{ with }\ g_{10000}(y)=y\, , \end{equation} 
and the rest of the orbit apparently develops in a tree-like manner. Theorem \ref{nonex} applies to $z$ and $y,$ and $\beta$ is determined by \eqref{last} as well as by the cycles of $\frac12 .$ However, the (upper) local dimensions of $\nu$ at $z$ and at $x=g_0(y)$ must coincide although the value of $x$ in the histogram is 1.4. The difference comes from the nonexponential growth given by the two cycles in \eqref{last} and determines a large factor $C$ in the ansatz $\nu(U(z,r))\approx C r^d$ but not a large dimension $d.$ The situation for $s'$ and $s"$ in Example 3 could be similar.\medskip

A recent result of Hare, Hare and Matthews \cite[Theorem 5.1]{HHM} says that for Pisot parameters local dimensions of periodic orbits are dense in the set of all possible local dimensions of $\nu .$
Our figures support the conjecture that supercritical intersections of periodic address curves are present at all parameters $t$ for which $\nu_t$ does not admit a bounded density. 
Although the last example indicates a difficulty, it has become clear that periodic orbits are important landmarks.  Of course there are other phenomena as invariant Cantor sets and multifractal spectra which were thoroughly studied for Pisot and Salem parameters \cite{LP, F4, F11, JSS, HHM}. Most promising in a study of the two-variable function $\Phi$ seems parameter dependence which we did not touch here. 
To mention just one aspect: points with unique address outside $\bf D$ have local dimension $1+\frac{\log 2}{\log\beta}$ \cite{ban}. It would be interesting to know two-dimensional local dimensions at intersection points of address curves.

\noindent
Christoph Bandt\\
Institute of Mathematics and Computer Science\\
University of Greifswald, Germany\\
\url{bandt@uni-greifswald.de}

\end{document}